\documentclass[12pt]{article}
\usepackage{microtype} \DisableLigatures{encoding = *, family = * }
\usepackage{subeqn,multirow}
\usepackage[pdftex]{graphicx}
\usepackage{amsfonts,amssymb}
\usepackage{natbib}
\bibpunct{(}{)}{;}{a}{,}{,}
\usepackage[bookmarksnumbered=true, pdfauthor={Wen-Long Jin}]{hyperref}

\newcommand{\commentout}[1]{}

\newcommand{\ba}{\begin{array}}
        \newcommand{\ea}{\end{array}}
\newcommand{\bc}{\begin{center}}
        \newcommand{\ec}{\end{center}}
\newcommand{\bdm}{\begin{displaymath}}
        \newcommand{\edm}{\end{displaymath}}
\newcommand{\bds} {\begin{description}}
        \newcommand{\eds} {\end{description}}
\newcommand{\ben}{\begin{enumerate}}
        \newcommand{\een}{\end{enumerate}}
\newcommand{\beq}{\begin{equation}}
        \newcommand{\eeq}{\end{equation}}
\newcommand{\bfg} {\begin{figure}}
        \newcommand{\efg} {\end{figure}}
\newcommand{\bi} {\begin {itemize}}
        \newcommand{\ei} {\end {itemize}}
\newcommand{\bpp}{\begin{pspicture}}
        \newcommand{\epp}{\end{pspicture}}
\newcommand{\bqn}{\begin{eqnarray}} 
        \newcommand{\eqn}{\end{eqnarray}}
\newcommand{\bqs}{\begin{eqnarray*}}
        \newcommand{\eqs}{\end{eqnarray*}}
\newcommand{\bsq}{\begin{subequations}}
        \newcommand{\esq}{\end{subequations}}
\newcommand{\bsl} {\begin{slide}[8.8in,6.7in]}
        \newcommand{\esl} {\end{slide}}
\newcommand{\bss} {\begin{slide*}[9.3in,6.7in]}
        \newcommand{\ess} {\end{slide*}}
\newcommand{\btb} {\begin {table}[h]}
        \newcommand{\etb} {\end {table}}
\newcommand{\bvb}{\begin{verbatim}}
        \newcommand{\evb}{\end{verbatim}}

\newcommand{\m}{\mbox}
\newcommand {\der}[2] {{\frac {\m {d} {#1}} {\m{d} {#2}}}}
\newcommand {\pd}[2] {{\frac {\partial {#1}} {\partial {#2}}}}

\newcommand{\cas}[1]{{{\left \{ \ba #1 \ea \right. }}}

\newcommand{\reff}[1] {{{Figure \ref {#1}}}}
\newcommand{\refe}[1] {{{(\ref {#1})}}}
\newcommand{\reft}[1] {{{\textbf{Table} \ref {#1}}}}

\newtheorem{theorem}{Theorem}[section]
\newtheorem{lemma}[theorem]{Lemma}
\newtheorem{corollary}[theorem]{Corollary}
\newtheorem{conjecture}[theorem]{Conjecture}

\def\a {{\alpha}}

\def\pmb#1{\setbox0=\hbox{$#1$}%
   \kern-.025em\copy0\kern-\wd0
   \kern.05em\copy0\kern-\wd0
   \kern-.025em\raise.0433em\box0 }
\def\bfxi{\pmb \xi}

\def\eop{{\hfill $\blacksquare$}}

\def\dt     {{\Delta t}}

\def\eop{{\hfill $\blacksquare$}}

\def\la {{{\lambda}}}
\def\o {{\omega}}

\def\FF{{\mathcal{F}}}
\def\L{{\textit{\textbf{L}}}}
\def\H{{\textit{\textbf{H}}}}

\usepackage{cleveref}

\begin {document}
\title{Continuous formulations and analytical properties of the link transmission model} 
\author{Wen-Long Jin\footnote{Department of Civil and Environmental Engineering, California Institute for Telecommunications and Information Technology, Institute of Transportation Studies, 4000 Anteater Instruction and Research Bldg, University of California, Irvine, CA 92697-3600. Tel: 949-824-1672. Fax: 949-824-8385. Email: wjin@uci.edu. Corresponding author}}

\maketitle

\begin{abstract}
The link transmission model (LTM) has great potential for simulating traffic flow in large-scale networks since it is much more efficient and accurate than the Cell Transmission Model (CTM). However, there lack general continuous formulations of LTM, and there has been no systematic study on its analytical properties such as stationary states and stability of network traffic flow. In this study we attempt to fill the gaps. First we apply the Hopf-Lax formula to derive Newell's simplified kinematic wave model with given boundary cumulative flows and the triangular fundamental diagram. We then apply the Hopf-Lax formula to define link demand and supply functions, as well as link queue and vacancy functions, and present two continuous formulations of LTM, by incorporating boundary demands and supplies as well as invariant macroscopic junction models. With continuous LTM, we define and solve the stationary states in a road network. We also apply LTM to directly derive a Poincar\'e map to analyze the stability of stationary states in a diverge-merge network. Finally we present an example to show that LTM is not well-defined with non-invariant junction models. We can see that Newell's model and LTM complement each other and provide an alternative formulation of the network kinematic wave model. This study paves the way for further extensions, analyses, and applications of LTM in the future.
 
\end{abstract}

{\bf Keywords}: Hopf-Lax formula; Newell's model; link transmission model; link demand and supply; link queue and vacancy sizes; invariant junction model; stationary states; stability.

\section{Introduction}
An understanding of congestion patterns in a road network is critical for developing efficient control, management, planning, and design strategies to improve safety, mobility, and environmental impacts of the transportation system. Network traffic flow models have been proposed at various scales: microscopic models \citep{hidas2005modelling}, the LWR model \citep{lighthill1955lwr,richards1956lwr}, the link queue model \citep{jin2012_link}, the neighborhood model \citep{daganzo2008analytical}, and continuum network models \citep{beckmann1952transportation,ho2006continuum}.

In particular, the celebrated LWR model, which can capture the initialization, propagation, and dissipation of vehicular queues on a road link through shock and rarefaction waves, has been extended as network kinematic wave models \citep{holden1995unidirection,garavello2006tfn}, by incorporating macroscopic merging, diverging, and junction models, which were first introduced in the Cell Transmission Model (CTM) and then established as analytical entropy conditions \citep{daganzo1995ctm,lebacque1996godunov,jin2012network}. 
CTM is a discrete Godunov version of the hyperbolic conservation law formulation of the network kinematic wave model, in which a link is divided into cells, a time duration discretized into time steps, boundary fluxes calculated from upstream demands and downstream supplies according to macroscopic junction models, and densities updated from the conservation law. Thus the computational cost is inversely proportional to the square of the time-step size \citep{leveque2002fvm}.
In contrast, \citet{yperman2006mcl,yperman2007link} introduced the link transmission model (LTM), which is another discrete version of the network kinematic wave model. In LTM, macroscopic junction models are also used to determine boundary fluxes, but the demand and supply functions are defined from cumulative flows based on Newell's formulation of the LWR model \citep{newell1993sim}; Newell's model was variational principle solutions to the Hamilton-Jacobian formulation of the LWR model in \citep{daganzo2005variationalKW,daganzo2005variationalKW2,daganzo2006variational}. But the computational cost of LTM is just inversely proportional to the time-step size; thus it is much more efficient than CTM. 
In addition, LTM was also shown to yield much more accurate solutions than CTM with the same time-step size.
Therefore, LTM has great potential in solving the network loading problem and the dynamical traffic assignment problem \citep{gentile2007macroscopic}.

However, in the discrete form, LTM has been mostly used as a simulation tool; and in the original LTM, all links are assumed to be initially empty. That is, there lack general continuous formulations of LTM, and there has been no systematic study on its analytical properties such as stationary states and stability of network traffic flow.

 In this study we attempt to fill the gaps: we will first derive two continuous formulations of LTM under arbitrary initial conditions and then apply them to solve stationary states and analyze their stability property in a road network. To derive the continuous formulations, we first apply the Hopf-Lax formula to derive Newell's model under general initial conditions \citep{evans1998pde,claudel2010lax2}, and then define link demand and supply functions as well link queue and vacancy sizes. With the two continuous formulations of LTM, we define stationary states and solve them for a road network with constant boundary conditions. Further, we analytically demonstrate that LTM can be unstable on a diverge-merge network. We also show that non-invariant junction models lead to ill-defined LTM. 

Note that, in \citep{han2012continuous}, the Hopf-Lax formula was used to derive a continuous formulation of LTM, but all links were assumed to be initially empty. In contrast, in this study we systematically apply the Hopf-Lax formula to derive Newell's model and two continuous formulations of LTM under general initial conditions. In addition, the junction models are more general here, and we also apply the continuous LTM to obtain analytical properties related to stationary states and their stability. Thus this study is more complete and leads to a better understanding of the properties of LTM.

The rest of the paper is organized as follows. In Section 2, we apply the Hopf-Lax formulation to derive Newell's model with a triangular fundamental diagram. In Section 3, we define link demand and supply functions from the Hopf-Lax formula and derive two continuous versions of LTM. In Section 4, we define stationary states and solve them in a network. We also show that non-invariant junction models lead to ill-defined LTM. In Section 5, we analytically demonstrate that LTM can be unstable on a diverge-merge network. Finally in Section 6, we discuss future research topics.

\section{The Hopf-Lax formula and Newell's model}
In the LWR model, the evolution of density, $k(x,t)$, speed, $v(x,t)$, and flux, $q(x,t)$, at location $x$ and time $t$ on a homogeneous road can be described based on the flow conservation equation, $k_t+q_x=0$, and a fundamental diagram $q=Q(k)$. Traditionally the LWR model has been written as the following hyperbolic conservation law:
\bqn
\pd{k}t+\pd{Q(k)}x&=&0, \label{lwr_hcl}
\eqn
where density is the state variable.

In \citep{newell1993sim}, the Hamilton-Jacobi formulation was introduced for the LWR model with a new state variable, $A(x,t)$, which is the cumulative flow passing $x$ before $t$ inside a spatial-temporal domain $\Omega$ and also known as a Moskowitz function \citep{moskowitz1965discussion}.
 Then $k=-A_x$, $q=A_t$, and the flow conservation equation is automatically satisfied when $A_{xt}=A_{tx}$. Then the fundamental diagram and, therefore, the LWR model, is equivalent to the following Hamilton-Jacobi equation
\bqn
A_t-Q(-A_x)&=&0. \label{lwr_vt}
\eqn 
Here the Hamiltonian is $\H(A_x)=-Q(-A_x)$, which is convex when the fundamental diagram is concave.

\subsection{The Hopf-Lax formula}
Many methods have been developed to uniquely solve the Hamilton-Jacobi equation when initial and boundary conditions in $A(x,t)$ are given.
In \citep{newell1993sim}, a minimization principle was used to solve \refe{lwr_vt} on a road link or network when on-ramp cumulative flows and off-ramp turning proportions are given. 
In \citep{daganzo2005variationalKW,daganzo2005variationalKW2}, the variational principle was applied to solve \refe{lwr_vt} analytically and numerically. 
Other solution methods, such as the Hopf-Lax formula, optimal control solutions, and viscous solutions, can be found in \citep{evans1998pde}. 
In this study, we will apply the Hopf-Lax formula to solve \refe{lwr_vt} under different types of boundary conditions: we first derive Newell's model when boundary cumulative flows are given, and then LTM when boundary demands and supplies are given. 

If we denote the Legendre transformation of $Q(k)$ by
\bqs
\L(u)&=&\sup_{k\in[0,K]} Q(k)-u \cdot k,
\eqs
where $K$ is the jam density, then \refe{lwr_vt} can be solved by the Hopf-Lax formula \citep[][Chapter 3]{evans1998pde}
\bqn
A(x,t)&=&\min_{(y,s)\in \partial \Omega(x,t)} A(y,s)+ (t-s) \L(\frac{x-y}{t-s} ), \label{hopf-lax}
\eqn
where $\partial \Omega$ is the boundary of $\Omega$, and $\partial \Omega(x,t)=\{(y,s)|(y,s)\in \partial \Omega, u=\frac{x-y}{t-s} \in[Q_k(K), Q_k(0)], t> s\}$ is the subset of boundary points that contribute to $A(x,t)$. 

In this study we only consider the triangular fundamental diagram as in \citep{newell1993sim}, 
\bqn
Q(k)&=&\min\{V k,  (K-k)W\},
\eqn
where $V$ is the free-flow speed, and $-W$ the shock wave speed in congested traffic. Thus the critical density $\bar{K}=\frac{W}{V+W}K$, and the capacity is $C=V \bar{K}$. Then the Lagrangian function is
\bqs
\L(u)&=&\sup_{k\in [0,K]} Q(k)-u\cdot k =(V-u) \bar{K}=C-\bar K u,
\eqs
for $u\in[-W,V]$.
The Hopf-Lax formula, \refe{hopf-lax}, can be written as
\bsq
\bqn
A(x,t)&=&\min_{(y,s)\in \partial \Omega(x,t)}  B(y,s;x,t), \label{hl-triangular}
\eqn
where $\partial \Omega(x,t)=\{(y,s)|(y,s)\in \partial \Omega, u=\frac{x-y}{t-s} \in[-W, V], t> s\}$, and  
\bqn
B(y,s;x,t)=A(y,s)+(t-s)\L(\frac{x-y}{t-s})=A(y,s)+  (t-s)C-(x-y)\bar{K}. \label{def:B}
\eqn
\esq 

We have the following special cases for $B(y,s;x,t)$: 
\ben
\item When $u=0$, $\L(0)=C$, and $B(y,s;x,t)=A(y,s)+(t-s)C$. 
\item When $u=V$, $\L(V)=0$, and $B(y,s;x,t)=A(y,s)$.
\item When $u=-W$, $\L(-W)=WK$, and $B(y,s;x,t)=A(y,s) + (t-s)KW=A(y,s)- (x-y)K$. 
\een
Furthermore, for any valid $u\in[-W,V]$ we can decompose the vector $(y,s)\to (x,t)$ ($t>s$) into two vectors with an intermediate point $(x_1,t_1)$ ($s\leq t_1\leq t$), such that $\frac{x_1-y}{t_1-s}=V$, $\frac{x-x_1}{t-t_1}=-W$, and $B(y,s;x,t)=A(y,s)-(x-x_1)K$.

Assume that $A(y,s)$ is given along a road segment between $x_1$ and $x_2>x_1$ at time $s$; i.e., $A(y,s)$ is known along $\partial \Omega_1=\{(y,s)|x_1\leq y\leq x_2\}$. Then for $(x,t)$ inside a cone defined by $t-s\geq \frac{x_2-x_1}{V+W}$ and $x_2-(t-s)W\leq x\leq x_1+(t-s)V$, $\partial \Omega_1$ is a subset of $\partial\Omega(x,t)$.  

\begin{lemma}  \label{x-lemma}
(i) If initially traffic is uncongested; i.e., if $k(y,s)=-\pd{A(y,s)}{y}\leq \bar{K}$ for $y\in[x_1,x_2]$, then 
\bqn
\min_{(y,s)\in \partial \Omega_1}  B(y,s;x,t)=B(x_1,s;x,t), \label{forward-wave}
\eqn
which is determined by the upstream end point;
(ii) If initially traffic is congested; i.e., if $k(y,s)> \bar{K}$ for $y\in[x_1,x_2]$, then
\bqn
\min_{(y,s)\in \partial \Omega_1}  B(y,s;x,t)=B(x_2,s;x,t), \label{backward-wave}
\eqn
which is determined by the downstream end point.
\end{lemma}
{\em Proof}. From \refe{def:B}, we have $\pd{B(y,s;x,t)}y=\pd{A(y,s)}{y}+\bar{K}$. Thus $B(y,s;x,t)$ increases in $y$, and \refe{forward-wave} is true, when $-\pd{A(y,s)}{y}\leq \bar{K}$. Thus $B(y,s;x,t)$ decreases in $y$, and \refe{backward-wave} is true, when $-\pd{A(y,s)}{y}> \bar{K}$. 
\eop

From Lemma \ref{x-lemma} we have the following corollary.
\begin{corollary} \label{x-corollary}
If the initial density satisfies
\bqn
k(y,s)\leq \bar{K}, \m{for } y\in [x_1,x_3] \qquad \m{and} \qquad k(y,s)> \bar{K}, \m{for } y\in (x_3,x_2] \label{regulartraffic}
\eqn
where $x_1\leq x_3\leq x_2$; i.e., if the upstream section is uncongested, and the downstream section congested,
then 
\bqn
\min_{(y,s)\in \partial \Omega_1} B(y,s;x,t)=\min \left\{B(x_1,s;x,t),B(x_2,s;x,t) \right\},
\eqn
which is determined by the two end points.
\end{corollary}
Note that the initial condition in \refe{regulartraffic} can be quite broad: traffic conditions can be uncongested ($x_3=x_2$) or congested ($x_3=x_1$) on the whole road segment, or uncongested in the upstream section and congested in the downstream section. From the traditional kinematic wave theory, we can see that rarefaction and shock waves can arise from the initial conditions; but \refe{regulartraffic} excludes congested upstream and uncongested downstream initial conditions, when transonic rarefaction waves can arise \citep{lebacque1996godunov}. 

Assume that $A(y,s)$ is given at $y$ during a time interval $[t_1,t_2]$; i.e., $A(y,s)$ is known along $\partial \Omega_2=\{(y,s)|t_1\leq s\leq t_2\}$.
Then for $(x,t)$ inside a cone defined by $t> t_2>t_1$ and $y-(t-t_2)W\leq x\leq y+(t-t_2)V$,  $\partial \Omega_2$ is a subset of $\partial\Omega(x,t)$. 

\begin{lemma} \label{t-lemma}
When $A(y,s)$ is given at $y$ between $t_1$ and $t_2>t_1$, we have 
\bqn
\min_{(y,s)\in \partial\Omega_2} B(y,s;x,t)=B(y,t_2;x,t). \label{bc}
\eqn
That is, the earlier boundary data is irrelevant.
\end{lemma}
{\em Proof}. From \refe{def:B}, we have $\pd{B(y,s;x,t)}s=\pd{A(y,s)}{s}-C\leq 0$. Thus $B(y,s;x,t)$ decreases in $s$, and \refe{bc} is true. 
\eop

\subsection{Newell's simplified kinematic wave model}
In \citep{newell1993sim}, Newell's simplified kinematic wave model was proposed to determine traffic conditions inside a homogeneous road section with a triangular fundamental diagram from initial and boundary cumulative flows.  That is, the Hamilton-Jacobi formulation of the LWR model, \refe{lwr_vt}, is solved with a minimum principle inside a U-shaped spatial-temporal domain $\Omega=[0,L]\times [0,\infty)$ as shown in \reff{u-domain}: where initial conditions, $N(x)=A(x,0)$ ($x\in [0,L]$), and boundary conditions, $F(t)=A(0,t)$ and $G(t)=A(L,t)$ ($t\geq 0$), are given. Equivalently, the initial density $k(x,0)$, in-flux $f(t)=\dot F(t)$, and out-flux $g(t)=\dot G(t)$ are given. 
The minimum principle was shown to be consistent with the variational principle in \cite{daganzo2005variationalKW} and can be used to solve \refe{lwr_vt} on inhomogeneous roads or a road network inside other spatial-temporal domains. However, all boundary cumulative flows have to be given.
That is, Newell's model can only handle boundary conditions of Dirichlet type.

\bfg\bc
\includegraphics[width=4in]{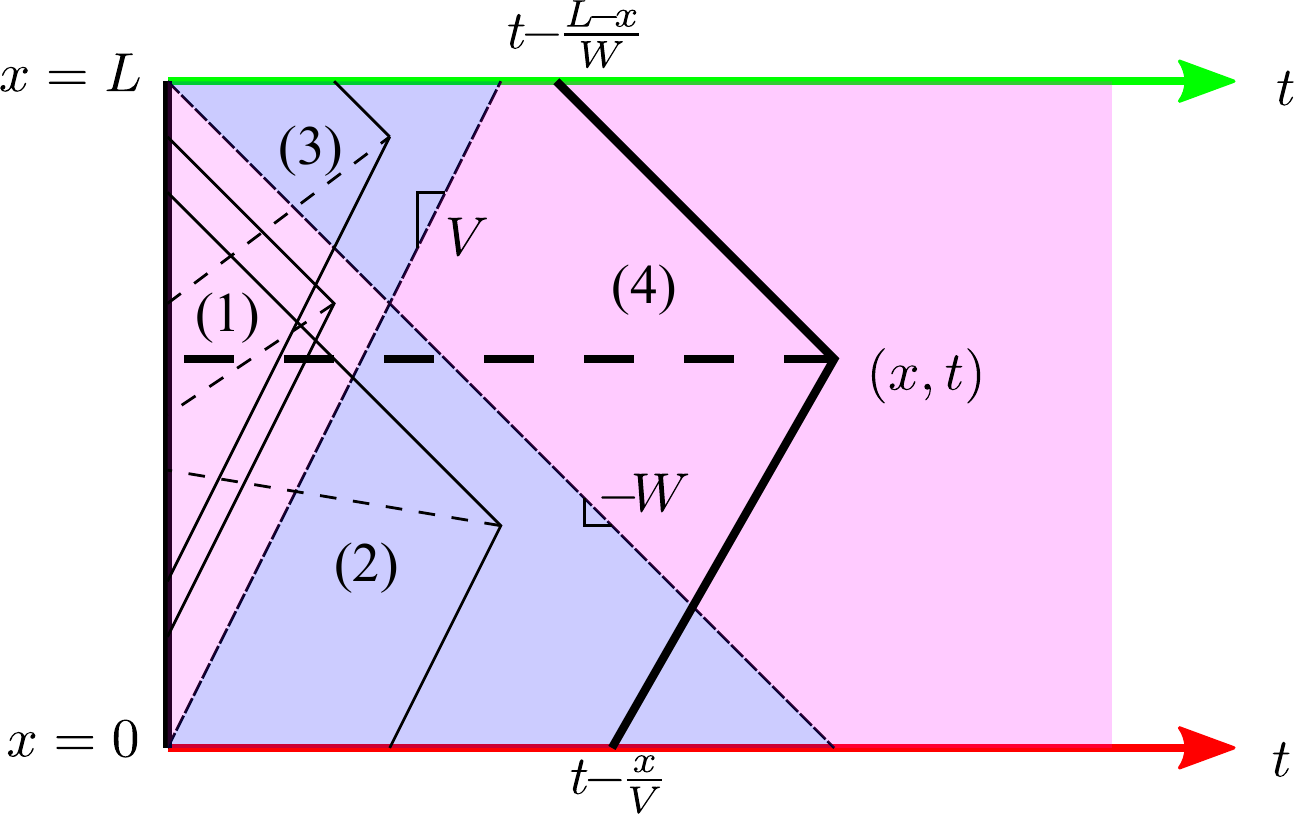}\caption{Newell's model for a U-shaped spatial-temporal domain with Dirichlet boundary conditions}\label{u-domain}
\ec\efg

 In the following, we apply the Hopf-Lax formula, \refe{hl-triangular}, to derive Newell's model and determine traffic conditions inside the U-shaped spatial-temporal domain, i.e., $A(x,t)$ for $x\in(0,L)$ and $t>0$.
For the four regions inside the domain, as shown in \reff{u-domain}, we have the following solutions.
\bsq \label{u-domain-sol}
\ben
\item For any point $(x,t)$ in region 1, we denote $x_1=x-V t$ and $x_2=x+W t$. Then \refe{hl-triangular} is equivalent to
\bqn
A(x,t)&=&\min_{y\in[x_1,x_2]} \{N(y)+Ct-(x-y)\bar{K}\}.
\eqn
\item For any point $(x,t)$ in region 2, we denote $x_2=x+W t$. From Lemma \ref{t-lemma}, \refe{hl-triangular} is equivalent to
\bqn
A(x,t)&=&\min_{y\in(0,x_2]} \{N(y)+Ct-(x-y)\bar{K},F(t-\frac{x}{V})\}.
\eqn
\item For any point $(x,t)$ in region 3, we denote $x_1=x-V t$. From Lemma \ref{t-lemma}, \refe{hl-triangular} is equivalent to
\bqn
A(x,t)&=&\min_{y\in[x_1,L)} \{N(y)+Ct-(x-y)\bar{K},G(t-\frac{L-x}{W})+ (L-x)K\}.
\eqn
\item For any point $(x,t)$ in region 4, from Lemma \ref{t-lemma}, \refe{hl-triangular} is equivalent to
\bqn
A(x,t)&=&\min_{y\in(0,L)} \{N(y)+Ct-(x-y)\bar{K},F(t-\frac{x}{V}),\nonumber \\&&G(t-\frac{L-x}{W})+ (L-x)K\}.
\eqn
\een
\esq
The solutions of $A(x,t)$ in Regions 2, 3, and 4 imply that
\bsq \label{feasiblecf}
\bqn
\m{Region 2: }F(t)&\leq& N(W t)+K W t, \\
\m{Region 3: }G(t)&\leq &N(L-Vt), \\
\m{Region 4: }F(t)&\leq & G(t-\frac{L}{W})+K L, \m{ and } G(t) \leq F(t- \frac{L}{V}),
\eqn
\esq
which are the necessary conditions for the initial-boundary problem to be well-posed.
That is, if \refe{feasiblecf} is violated, then Newell's simplified model is not well defined. This is consistent with the observation that the Cauchy-Dirichlet problem may not be well-posed under general initial and boundary conditions for a hyperbolic conservation law \citep{bardos1979first}.

In the following we present some special solutions to Newell's model.
\ben 
\item When \refe{regulartraffic} is satisfied for the initial traffic density; i.e., when there is no transonic rarefaction wave initially, from Corollary \ref{x-corollary}  we have the following simplified version of \refe{u-domain-sol}:
\bsq \label{u-domain-sol1}
\bqn
\m{Region 1: }A(x,t)&=&\min\{N(x-Vt),N(x+W t)+K W t\}, \\
\m{Region 2: }A(x,t)&=&\min\{F(t-\frac{x}{V}),N(x+W t)+K W t\}, \\
\m{Region 3: }A(x,t)&=&\min\{N(x-Vt),G(t-\frac{L-x}{W})+ (L-x)K\}, \\
\m{Region 4: }A(x,t)&=&\min\{F(t-\frac{x}{V}),G(t-\frac{L-x}{W})+ (L-x)K\}. \label{newell-original}
\eqn
\esq
In some studies \citep[e.g.][]{laval2012stochastic},  \refe{newell-original} in region 4 is referred to as Newell's simplified kinematic wave model. Note that the solution of $A(x,t)$ in region 4 depends on the initial conditions when \refe{regulartraffic} is violated.

\item When the initial traffic density is constant at $k_0$, then $F(0)=G(0)+k_0L$, and $N(x)=G(0)+k_0(L-x)$. Thus we have the following solutions:
\bsq \label{u-domain-sol-constant}
\bqn
\m{Region 1: }A(x,t)&=&F(0)-k_0x+\min\{k_0Vt,(K-k_0) W t\}, \\
\m{Region 2: }A(x,t)&=&\min\{F(t-\frac{x}{V}),F(0)-k_0x+(K-k_0) W t\}, \\
\m{Region 3: }A(x,t)&=&\min\{F(0)+(Vt-x)k_0,G(t-\frac{L-x}{W})+ (L-x)K\}, \\
\m{Region 4: }A(x,t)&=&\min\{F(t-\frac{x}{V}),G(t-\frac{L-x}{W})+ (L-x)K\}. 
\eqn
\esq

\item When the road is initially empty, then $N(x)=F(0)=G(0)$, and \refe{u-domain-sol1} is equivalent to
\bsq \label{u-domain-sol2}
\bqn
\m{Region 1: }A(x,t)&=&G(0), \\
\m{Region 2: }A(x,t)&=&F(t-\frac{x}{V}), \\
\m{Region 3: }A(x,t)&=&G(0), \\
\m{Region 4: }A(x,t)&=&\min\{F(t-\frac{x}{V}),G(t-\frac{L-x}{W})+K (L-x)\}.
\eqn
\esq

\een

\section{Continuous formulations of the link transmission model}

\btb
\bc
\begin{tabular}{|c|l||c|l|}\hline
\multicolumn{4}{|c|}{Sets}\\\hline
$I$ & set of origins & $O$ & set of destinations\\\hline
$A$ & set of regular links & $J$ & set of junctions\\\hline
$A'$ & set of all links & $ P $ & set of commodities in the network \\\hline
$I_j$ & set of upstream links of junction $j$ & $O_j$ & set of downstream links of junction $j$ \\\hline
$ P _a$ & set of commodities using link $a\in A'$ & $\Delta$& unidirectional road network \\\hline \hline
\multicolumn{4}{|c|}{Constants for a homogeneous regular link $a$}\\\hline
$L_a$& length  & $x_a$ & location \\\hline
$V_a$& free-flow speed & $-W_a$ & shock wave speed in congested traffic \\\hline
$K_a$& jam density & $C_a$ & capacity \\\hline\hline
\multicolumn{4}{|c|}{Variables for a homogeneous regular link $a$}\\\hline
$k_a(x_a,t)$&density &$q_a(x_a,t)$&flux \\\hline
$A_a(x_a,t)$& cumulative flow & $N_a(x_a)$ & initial cumulative flow  \\\hline
$F_a(t)$ & cumulative in-flow & $G_a(t)$ &cumulative out-flow \\\hline
$f_a(t)$ & in-flux & $g_a(t)$ &out-flux\\\hline
$d_a(t)$ & link demand & $s_a(t)$ & link supply \\\hline
$\la_a(t)$ & link queue size & $\gamma_a(t)$ & link vacancy size \\\hline\hline
\multicolumn{4}{|c|}{Variables for commodity $ p $ on a homogeneous regular link $a$}\\\hline
$F_{a, p }(t)$ & cumulative in-flow & $G_{a, p }(t)$ & cumulative out-flow  \\\hline
$f_{a, p }(t)$ & in-flux  & $g_{a, p }(t)$ & out-flux  \\\hline
$\eta_{a, p }(t)$ & upstream proportion  &$\xi_{a, p }(t)$ & downstream proportion  \\\hline\hline
\multicolumn{4}{|c|}{Variables for a junction $j$}\\\hline
$\theta_j(t)$ & critical demand level & $\xi_{a\to b}(t)$ & turning proportion \\\hline
\end{tabular}
\ec
\caption{A list of notations for the multi-commodity link transmission model}\label{listofnotations}
\etb

In \citep{yperman2006mcl,yperman2007link}, a discrete LTM was introduced as follows: inspired by Newell's simplified kinematic wave model, they first defined link demand (sending flow) and supply (receiving flow) in boundary cumulative flows  
and then incorporated them into macroscopic junction models, which were first proposed in the Cell Transmission Model (CTM) \citep{daganzo1995ctm,lebacque1996godunov}. 
Note that LTM is different from CTM, since demand and supply functions are defined in cell densities in CTM but boundary cumulative flows in LTM; it is also different from Newell's model, since Newell's model can only handle Dirichlet boundary conditions but LTM can handle  the Bardos-Leroux-Nedelec (BLN) and periodic boundary conditions through demand and supply functions and junction models as in CTM \citep{bardos1979first,lebacque2005network}. In this sense, LTM is a network extension of Newell's model and an alternative formulation of the network kinematic wave model.

\begin{figure}\bc
\includegraphics[width=3in]{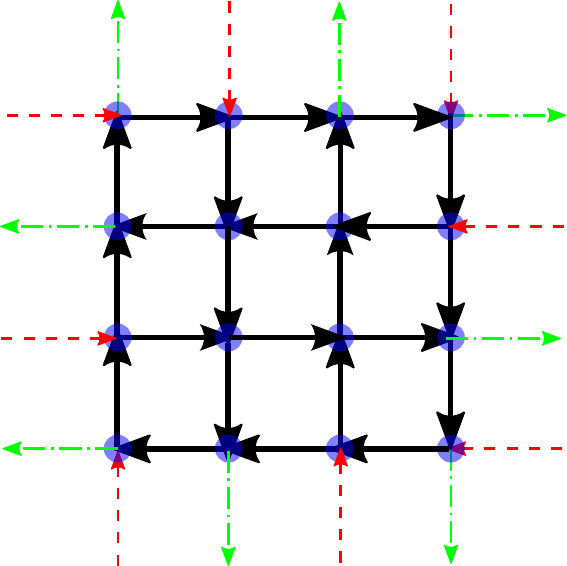}\caption{A $4\times 4$ grid network}\label{gridnetwork}
\ec
\end{figure}

However, no continuous formulations of link demand and supply functions have been rigorously defined under general initial conditions. In this section, we will apply the Hopf-Lax formula to derive two continuous versions of LTM under general initial conditions for a general road network, e.g., a grid network shown in \reff{gridnetwork}. 
In a network, we denote the sets of origins (dashed red lines), destinations (dash-dotted green lines), regular links (solid black lines), and junctions (blue dots) by $I$, $O$, $A$, and $J$, respectively. Here origins and destinations can be treated as dummy links of zero length, and we denote $A'=I\cup O\cup A$.
We assume that movements from all upstream links to all downstream links are permitted. We denote the set of the upstream links of junction $j$ by $I_j$ and the set of the downstream links by $O_j$. 

In addition, vehicles are categorized into commodities according to their paths; i.e., vehicles on the same path belong to the same commodity. \footnote{In emergency evacuation situations, vehicles who comply with evacuation orders may not have pre-defined routes, and the commodity of such vehicles can use many paths; if multi-class vehicles are considered, then vehicles of two commodities can share the same path.} 
The set of commodities in the whole network is denoted by $ P $, and the set of commodities using link $a\in A'$ is denoted by $ P _a$. 
Then a unidirectional road network  can be defined by
\bqn
\Delta&=&\left(I, O, A, J, \{( I_j, O_j): j\in J\} , \{ P _a:a\in A'\} \right)
\eqn

For a network $\Delta$, the constants and variables are defined in \reft{listofnotations}. The traffic dynamics on link $a$ are described by the LWR model, \refe{lwr_hcl} or \refe{lwr_vt}, with a triangular fundamental diagram:
\bqn
q_a(x_a,t)&=&Q_a(k_a(x_a,t))=\min\{V_a k_a (x_a,t), (K_a-k_a(x_a,t)) W_a\},
\eqn
where $x_a\in[0,L_a]$ is the coordinate for link $a$.

\subsection{Link demand and supply functions}
\begin{figure} \bc
$\ba{r@{\hspace{0.3in}}r}
\includegraphics[height=2in]{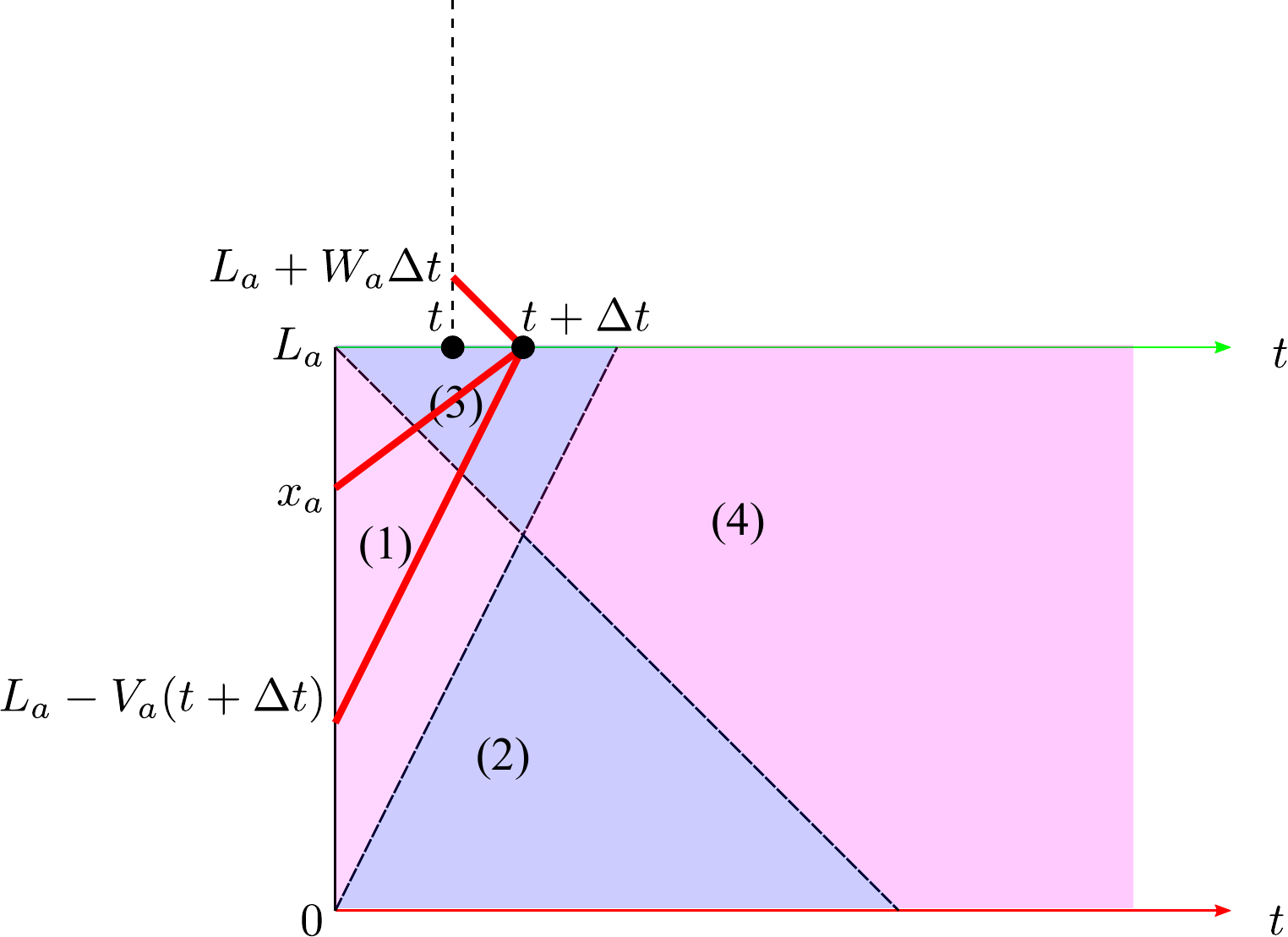} &
\includegraphics[height=2in]{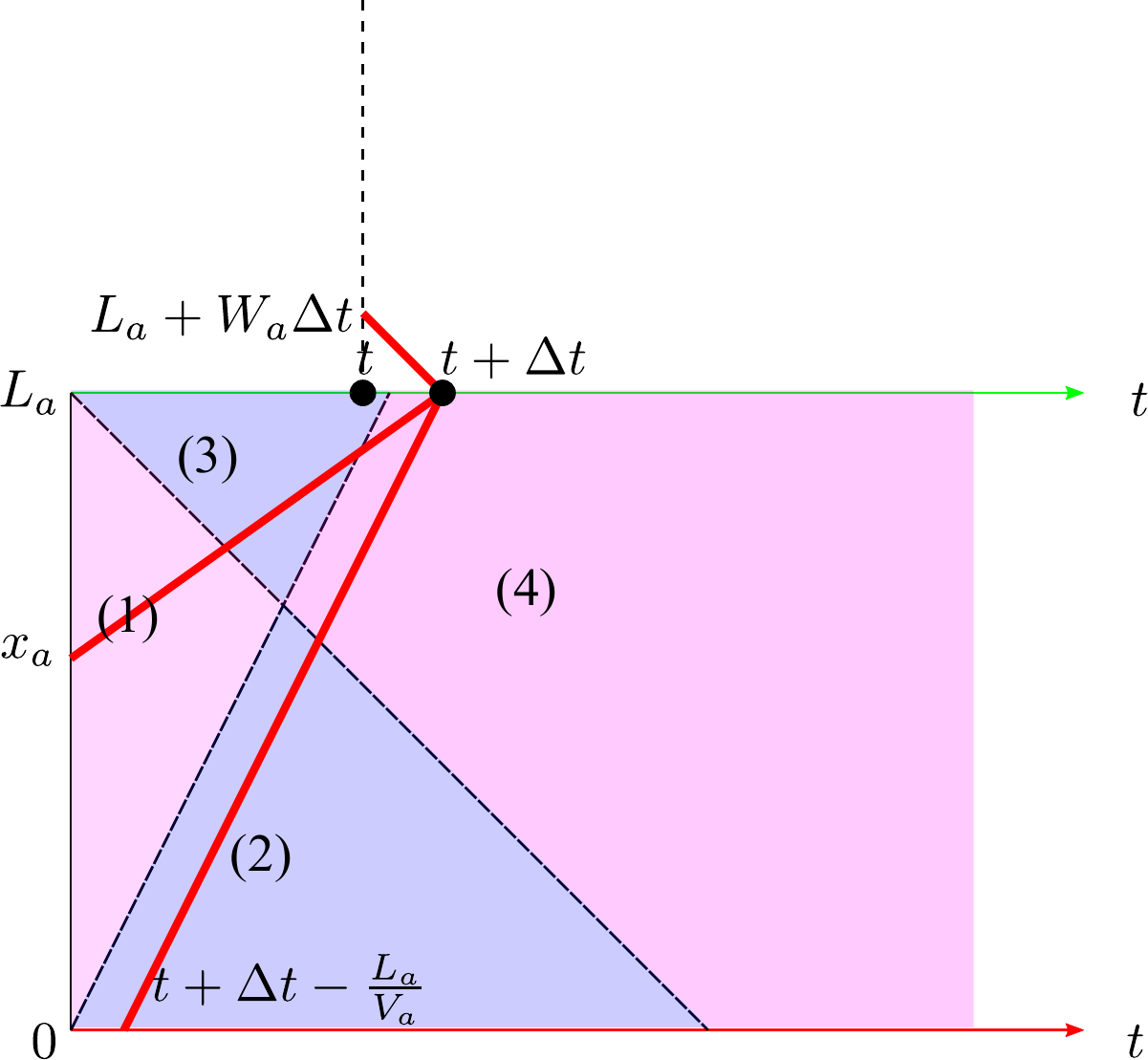} \\
\multicolumn{1}{c}{\mbox{\bf (a)}} &
    \multicolumn{1}{c}{\mbox{\bf (b)}}\\
\includegraphics[height=2in]{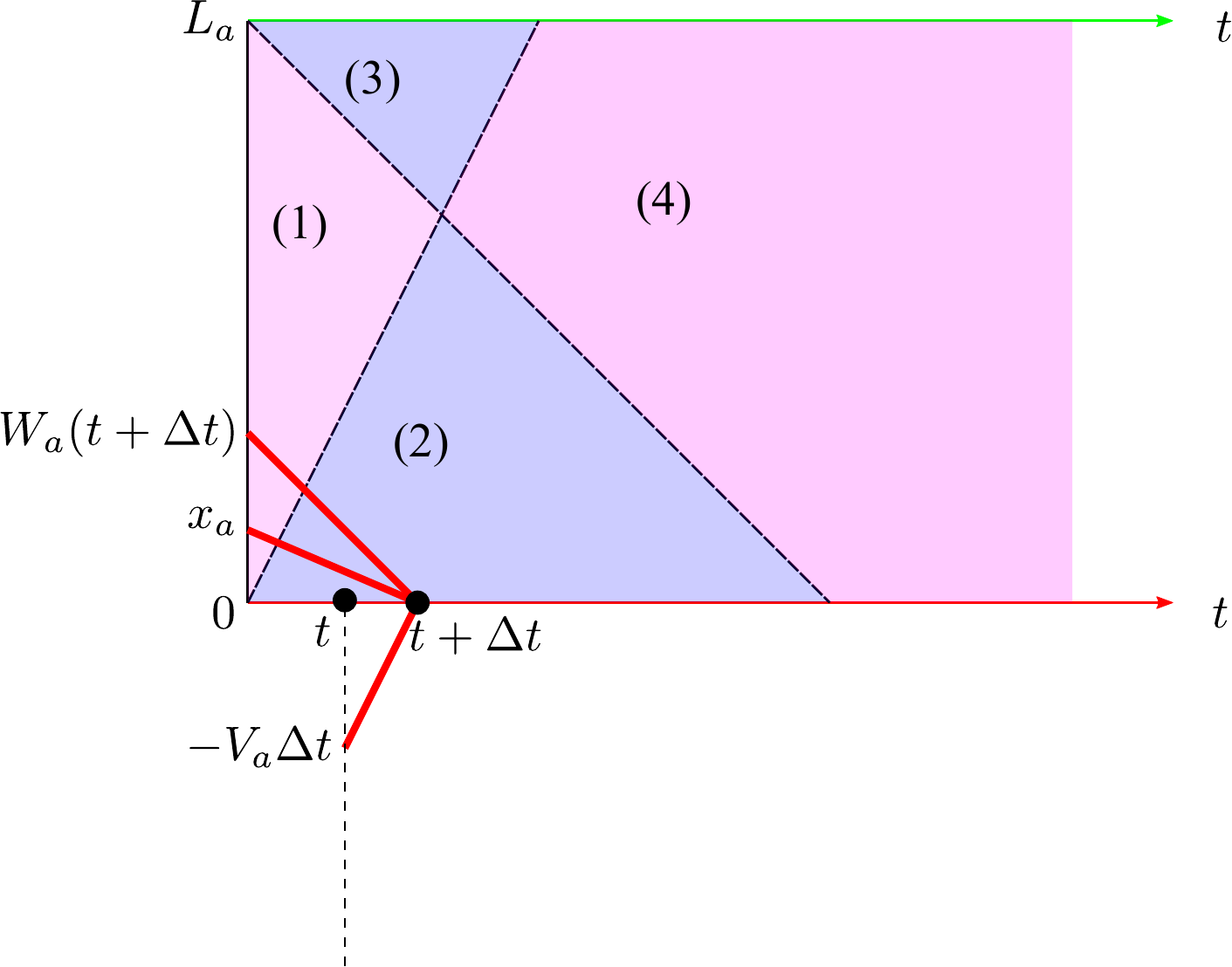} &
\includegraphics[height=2in]{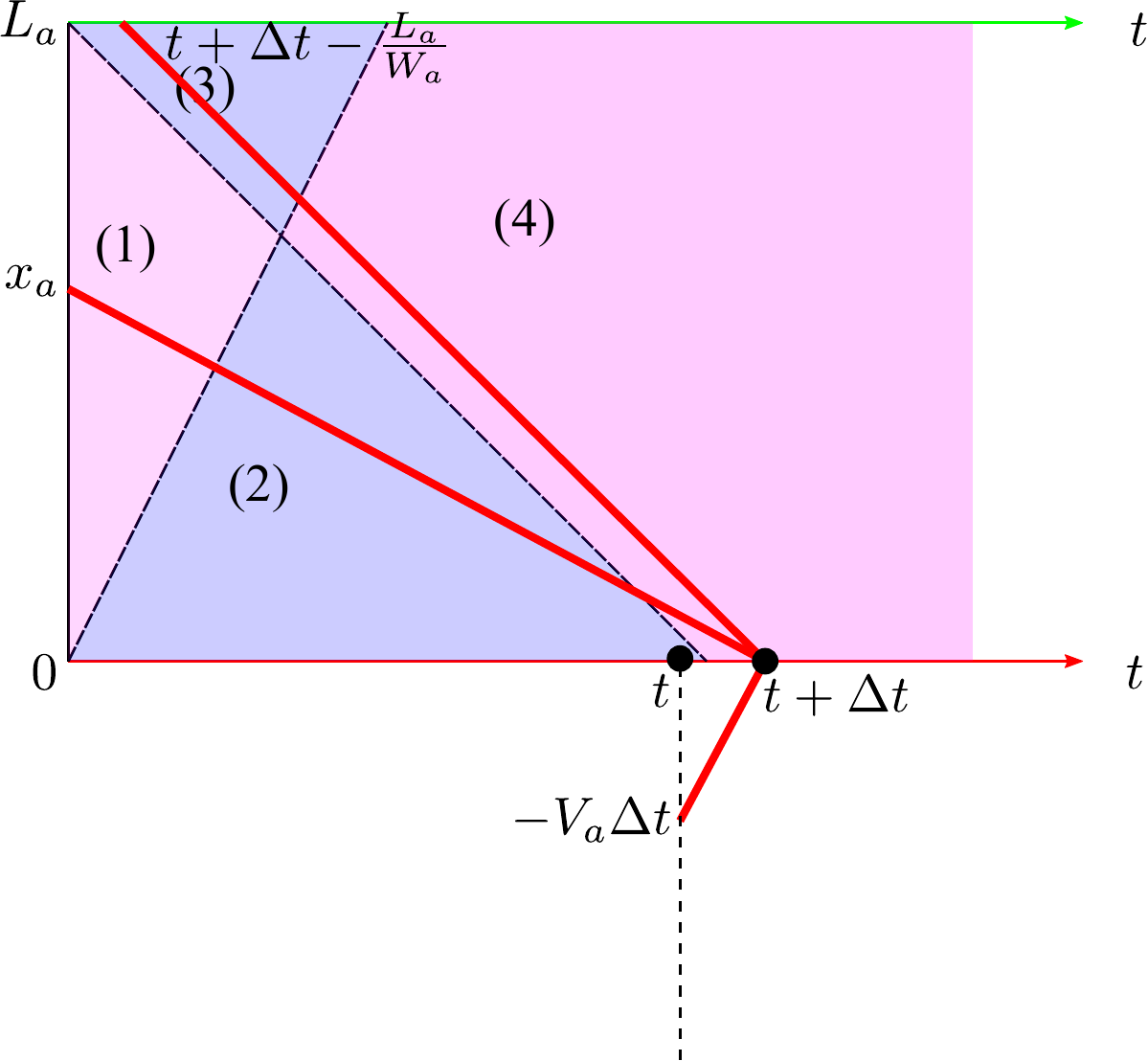} \\
\multicolumn{1}{c}{\mbox{\bf (c)}} &
    \multicolumn{1}{c}{\mbox{\bf (d)} }    
\ea$
\caption{Definitions of link demands and supplies} \label{fig:ds} \ec 
\end{figure}

We follow \citep{daganzo1995ctm,lebacque1996godunov} and define the link demand function at time $t$, $d_a(t)$, as the out-flux with an empty downstream link, and the link supply function at time $t$, $s_a(t)$, as the in-flux with a jammed upstream link. Here we assume that $F_a(\tau)$ and $G_a(\tau)$ are known for $\tau\in[0,t]$, and $N_a(x_a)$ are known for $x_a\in[0,L_a]$.
\ben
\item To define the link demand, we extend link $a$ such that $x_a\in[0,\infty)$, but the downstream link is empty; i.e., $k_a(x_a,t)=0$ for $x_a>L_a$. In this case, we first apply the Hopf-Lax formula to calculate the out-flow $\hat G_a(t+\dt)$ for a small $\dt>0$ with a pan-shaped spatial-temporal domain shown in \reff{fig:ds}(a) and (b), and then define the link demand $d_a(t)$ by
\bsq
\bqn
d_a(t)&=&\lim_{\dt\to 0^+}\frac{\hat G_a(t+\dt)-G_a(t)}{\dt}. \label{def:d}
\eqn
From Lemma \ref{t-lemma}, we can obtain $\hat G_a(t+\dt)$ as follows:
\ben
\item If $t+\dt\leq \frac{L_a}{V_a}$, as shown in \reff{fig:ds}(a), 
\bqn
\hat G_a(t+\dt)&=&\min_{x_a\in[L_a-(t+\dt)V_a,L_a)}\{B(x_a,0;L_a,t+\dt),G_a(t)+C_a \dt\}. 
\eqn

\item  If $t+\dt> \frac{L_a}{V_a}$, as shown in \reff{fig:ds}(b), 
\bqn
\hat G_a(t+\dt)&=&\min_{x_a\in(0,L_a)}\{F_a(t+\dt-\frac{L_a}{V_a}),\nonumber \\ &&B(x_a,0;L_a,t+\dt),G_a(t)+C_a \dt\}. 
\eqn
\een
\esq
As expected, $\hat G_a(t+\dt)$ only depends on traffic conditions inside link $a$. Note that $\hat G_a(t+\dt)$ is the ideal value of $G_a(t+\dt)$ when the downstream is empty, but the two values may not be the same in reality when the downstream is not.

\item  To define the link supply, we extend link $a$ such that $x_a\in(-\infty,0]$, but upstream link is jammed; i.e., $k_a(x_a,t)=K_a$ for $x_a<0$. In this case, we first apply the Hopf-Lax formula to calculate the in-flow $\hat F_a(t+\dt)$ for a small $\dt>0$ with a pan-shaped spatial-temporal domain shown in \reff{fig:ds}(c) and (d), and then define the link supply $s_a(t)$ by
\bsq
\bqn
s_a(t)&=&\lim_{\dt\to 0^+}\frac{\hat F_a(t+\dt)-F_a(t)}{\dt}. \label{def:s}
\eqn
From Lemmas \ref{x-lemma} and \ref{t-lemma}, we can obtain $\hat F_a(t+\dt)$ as follows: 
\ben
\item If $t+\dt\leq \frac{L_a}{W_a}$, as shown in \reff{fig:ds}(c), 
\bqn
\hat F_a(t+\dt)&=&\min_{x_a\in(0,(t+\dt)W_a]}\{B(x_a,0;0,t+\dt),F_a(t)+C_a \dt\}. 
\eqn

\item  If $t+\dt> \frac{L_a}{W_a}$, as shown in \reff{fig:ds}(d), 
\bqn
\hat F_a(t+\dt)&=&\min_{x_a\in(0,L_a)}\{G_a(t+\dt-\frac{L_a}{W_a})+K_a L_a,\nonumber \\ && B(x_a,0;0,t+\dt),F_a(t)+C_a \dt\}. 
\eqn
\een
\esq
As expected, $\hat F_a(t+\dt)$ only depends on traffic conditions inside link $a$. Note that $\hat F_a(t+\dt)$ is the ideal value of $F_a(t+\dt)$ when the upstream is jammed, but the two values may not be the same in reality when the upstream is not jammed.
\een

In the discussions hereafter, we assume that \refe{regulartraffic} is satisfied; i.e., there exists no initial transonic rarefaction wave on a link. Then from Corollary \ref{x-corollary} $\hat G_a(t+\dt)$ and $\hat F_a(t+\dt)$ can be simplified as follows: 
\bqs
\hat G_a(t+\dt)&=&\cas{{ll} \min\{ N_a(L_a-V_a(t+\dt)),G_a(t)+C_a \dt\}, & t+\dt \leq \frac{L_a}{V_a}\\
\min\{ F_a(t+\dt-\frac{L_a}{V_a}),G_a(t)+C_a \dt\}, & t+\dt > \frac{L_a}{V_a}}\\
\hat F_a(t+\dt)&=&\cas{{ll} \min\{ N_a(W_a(t+\dt))+ (t+\dt)K_a W_a,F_a(t)+C_a \dt\}, & t+\dt\leq \frac{L_a}{W_a}\\\min\{ G_a(t+\dt-\frac{L_a}{W_a})+K_a L_a,F_a(t)+C_a \dt\}, & t+\dt> \frac{L_a}{W_a}}
\eqs
Then from \refe{def:d} and \refe{def:s}, the link demand and supply can be re-written as:
\bsq \label{def:ds}
\bqn
d_a(t)&=&\cas{{ll} \min\left\{k_a(L_a-V_a t,0) V_a+H(\la_a(t)), C_a\right\}, & t\leq \frac {L_a}{V_a}\\
\min\left\{ f_a(t-\frac{L_a}{V_a}) +H(\la_a(t)), C_a \right\}, & t>\frac{L_a}{V_a}}\\
s_a(t)&=&\cas{{ll} \min\left\{K_a W_a -k_a(W_a t,0) W_a+H(\gamma_a(t)), C_a\right\}, & t\leq \frac{L_a}{W_a}\\  \min\left\{ g_a(t-\frac{L_a}{W_a}) +H(\gamma_a(t)), C_a \right\}, &t>\frac{L_a}{W_a}}
\eqn
\esq
where the indicator function $H(y)$ for $y\geq 0$ is defined as
\bqn
H(y)&=&\lim_{\dt\to 0^+}\frac{y}{\dt}=\cas{{ll}0, &y=0;\\+\infty, &y>0.} \label{indicatorfunction}
\eqn
the in- and out-fluxes for link $a$ are 
\bsq
\bqn 
f_a(t)&=&\der{}{t}F_a(t),\\
g_a(t)&=&\der{}{t}G_a(t),
\eqn
\esq
and two new variables, including the link queue size, $\la_a(t)$, and the link vacancy size, $\gamma_a(t)$, are defined as follows:
\bsq \label{queuevacancy}
\bqn
\la_a(t)&=& \cas{{ll} N_a(L_a-V_a t)-G_a(t), & t\leq \frac {L_a}{V_a} \\F_a(t-\frac{L_a}{V_a})-G_a(t), & t>\frac{L_a}{V_a}}\\
\gamma_a(t)&=& \cas{{ll} N_a(W_a t)+K_a W_a t-F_a(t), & t\leq \frac{L_a}{W_a}\\ G_a(t-\frac{L_a}{W_a})+K_a L_a-F_a(t), &t>\frac{L_a}{W_a}}
\eqn
\esq
From \refe{queuevacancy}, we have $\la_a(0)=\gamma_a(0)=0$; from \refe{feasiblecf}, we can see that both $\la_a(t)$ and $\gamma_a(t)$ are non-negative. 

In the following we consider two special cases of \refe{regulartraffic}: 
\ben
\item When the initial traffic density is constant on link $a$; i.e., when $k_a(x_a,0)=k_a$ for $x_a\in[0,L_a]$, we have
\bsq
\bqn
\la_a(t)&=& \cas{{ll} k_a V_a t+G_a(0) -G_a(t), & t\leq \frac {L_a}{V_a} \\F_a(t-\frac{L_a}{V_a})-G_a(t), & t>\frac{L_a}{V_a}}\\
\gamma_a(t)&=& \cas{{ll} (K_a-k_a)  W_a t+F_a(0)-F_a(t), & t\leq \frac{L_a}{W_a}\\ G_a(t-\frac{L_a}{W_a})+K_a L_a-F_a(t), &t>\frac{L_a}{W_a}}\\
d_a(t)&=&\cas{{ll} \min\left\{k_a V_a+H(\la_a(t)), C_a\right\}, & t\leq \frac {L_a}{V_a}\\
\min\left\{ f_a(t-\frac{L_a}{V_a}) +H(\la_a(t)), C_a \right\}, & t>\frac{L_a}{V_a}}\\
s_a(t)&=&\cas{{ll} \min\left\{(K_a-k_a) W_a+H(\gamma_a(t)), C_a\right\}, & t\leq \frac{L_a}{W_a}\\  \min\left\{ g_a(t-\frac{L_a}{W_a}) +H(\gamma_a(t)), C_a \right\}, &t>\frac{L_a}{W_a}}
\eqn
\esq

\item When road $a$ is initially empty; i.e., when $k_a(x_a,0)=0$ for $x_a\in[0,L_a]$, 
we have

$\la_a(t)=0$ for $t\leq \frac{L_a}{V_a}$, and $d_a(t)=0$ for $t\leq \frac{L_a}{V_a}$; $\gamma_a(t)\geq (K_aW_a-C_a)t\geq 0$ for $t\leq \frac{L_a}{W_a}$, and $s_a(t) \geq \min\{K_a W_a, C_a\}=C_a$ for $t\leq \frac{L_a}{W_a}$. 
In this case, from \refe{queuevacancy} and \refe{def:ds} we have the following demand and supply functions:
\bqs
\la_a(t)&=& \cas{{ll} 0, & t\leq \frac {L_a}{V_a} \\F_a(t-\frac{L_a}{V_a})-G_a(t), & t>\frac{L_a}{V_a}}\\
\gamma_a(t)&=& \cas{{ll} K_a  W_a t+F_a(0)-F_a(t), & t\leq \frac{L_a}{W_a}\\ G_a(t-\frac{L_a}{W_a})+K_a L_a-F_a(t), &t>\frac{L_a}{W_a}}\\
d_a(t)&=&\cas{{ll} 0,& t \leq \frac{L_a}{V_a}; \\ f_a(t-\frac{L_a}{V_a}), & \m{if } F_a(t-\frac{L_a}{V_a})=G_a(t), t > \frac{L_a}{V_a};\\ C_a, & \m{if } F_a(t-\frac{L_a}{V_a})>G_a(t), t > \frac{L_a}{V_a}; }\\
s_a(t)&=&\cas{{ll} C_a, & t \leq \frac{L_a}{W_a}; \\g_a(t-\frac{L_a}{W_a}), & \m{if } F_a(t)=G_a(t-\frac{L_a}{W_a})+K_aL_a,t > \frac{L_a}{W_a};\\C_a, &\m{if }F_a(t)< G_a(t-\frac{L_a}{W_a})+K_aL_a, t > \frac{L_a}{W_a}.}
\eqs
Here the demand and supply functions on an initially empty road are consistent with those in \citep{han2012continuous}. In addition, the demand function is the same as that in the spatial queue model \citep{nie2002sqm,zhang2013modelling}, but not the supply function. Thus the spatial queue model is not the same as LTM.
\een

\subsection{An invariant junction model}
In LTM, macroscopic junction models can be used to determine boundary fluxes at a general junction $j$ from the upstream demands $d_a(t)$ ($a\in I_j$),  downstream supplies $s_b(t)$ ($b\in O_j$), and  turning proportions $\xi_{a\to b}(t)$, where $\xi_{a\to b}(t)\in [0,1]$ and $\sum_{b\in O_j} \xi_{a\to b}(t)=1$ for any $a$. 
Here we denote a general junction model by
\bqn
({\bf g}_j(t),{\bf f}_j(t))&=&\FF({\bf d}_j(t),{\bf s}_j(t), \bfxi_j(t)), \label{junctionfluxes}
\eqn
where  ${\bf d}_j(t)$ is the set of upstream commodity demands,  ${\bf s}_j(t)$ the set of downstream supplies, $\bfxi_j(t)$ the matrix of turning proportions, ${\bf g}_j(t)$ the set of out-fluxes from all upstream links, and ${\bf f}_j(t)$ the set of in-fluxes to all downstream links. 
Physically, such a junction model is determined by the characteristics of the junction bottleneck as well as vehicles' merging and diverging behaviors; mathematically, it is the so-called entropy condition that is used to pick out the unique weak solution for a system of hyperbolic conservation laws \citep{jin2012network}.

From the definitions of link demand and supply functions as well as conservation of traffic at the junction, the boundary fluxes should satisfy the following conditions:
\bsq \label{necessaryjunction}
\bqn
0&\leq & g_a(t)\leq  d_a(t), \quad a\in I_j,\\
0&\leq &f_b(t)\leq  s_b(t), \quad b\in O_j,\\
\sum_{a\in I_j} g_a(t)&=&\sum_{b\in O_j} f_b(t).
\eqn
Further, if vehicles follow the first-in-first-out diverging principle \citep{daganzo1995ctm}, we then have
\bqn
f_b(t)&=&\sum_{a\in I_j} g_a(t) \xi_{a\to b}(t).
\eqn
\esq
However, \refe{necessaryjunction} is not sufficient to uniquely determine the boundary fluxes. For example, for a linear junction connecting link $1$ to link $2$, the boundary flux can be any number between $0$ and $\min\{d_1(t), s_2(t)\}$. 
Therefore, additional rules are needed. 
Here we assume that all vehicles follow the fair merging rule \citep{jin2010merge} and the total flux is maximized. 

Then we obtain the following junction model, which was first derived in \citep{jin2012kinematic}. 
First we define the critical demand level $\theta_j(t)$ by the following min-max problem
\bsq\label{junctionmodel}
\bqn
\theta_j(t)&=&\min\{1,\min_{b\in O_j}  \max_{I_1(t)} \frac{s_b(t)-\sum_{\a\in I_j\setminus I_1(t)}d_\a(t) \xi_{\a\to b}(t)}{\sum_{a\in I_1(t)} C_a \xi_{a\to b}(t)}\}. \label{criticaldl}
\eqn
Here $I_1(t)$ a non-empty subset of $I_j$.
Then we calculate the out-flux of upstream link $a\in  I_j$ by
\bqn
g_a(t)&=&\min\{d_a(t),\theta_j(t) C_a\}, \label{outflux}
\eqn
and the in-flux of downstream link $b\in  O_j$ by
\bqn
f_b(t)&=&\sum_{a\in  I_j} g_a(t) \xi_{a\to b}(t). \label{influx}
\eqn
In addition, the commodity-flux is 
\bqn
f_{b,\o}(t)=g_{a,\o}(t)=g_a(t) \xi_{a,\o}(t).
\eqn
\esq

In \citep{jin2012_riemann}, it was shown that this junction model is invariant \citep{lebacque2005network}.
In Section 6 we will show that non-invariant junction models cannot be used in LTM.
Note that general junction models were also obtained in other studies \citep{tampere2011generic,flotterod2011operational}. But no explicit formulas were derived to calculate boundary fluxes, or the models were not rigorously proved to be invariant. 

\subsection{Two continuous formulations of the link transmission model}
Combining the definitions of demand and supply functions, the definitions of queue and vacancy sizes, the junction model, and the relationship between cumulative flows and fluxes, we can then derive two continuous formulations of LTM.

In the first formulation of the continuous LTM, we use the cumulative flows, $F_a(t)$ and $G_a(t)$, as the unknown state variables. We have the following evolution equations:
\bsq\label{cltm-1}
\bqn
\der{F_a(t)}t&=&f_a(t),\quad t\geq 0\\
\der{G_a(t)}t&=&g_a(t), \quad t\geq 0
\eqn
\esq
where the link queue and vacancy sizes are given by \refe{queuevacancy}, the link demand and supply functions by \refe{def:ds}, and $f_a(t)$ and $g_a(t)$ can by the junction model \refe{junctionmodel}. 
In particular, when link $a$ is initially empty, the link demand and supply can be calculated from link cumulative flows at $t$ and earlier times:
\bqs
d_a(t)&=&\cas{{ll} 0, & t\leq \frac{L_a}{V_a}\\\min\{\der{F_a(t-\frac{L_a}{V_a})}t+H(F_a(t-\frac{L_a}{V_a})-G_a(t)), C_a\}, & t>\frac{L_a}{V_a}}\\
s_a(t)&=&\cas{{ll} C_a, &t\leq \frac{L_a}{W_a}\\\min\{\der{G_a(t-\frac{L_a}{W_a})}t+H(G_a(t-\frac{L_a}{W_a})+K_aL_a-F_a(t)), C_a\}, &  t>\frac{L_a}{W_a}}
\eqs

In the second formulation of the continuous LTM, we use the link queue and vacancy sizes, $\la_a(t)$ and $\gamma_a(t)$, as the unknown state variables. Differentiating \refe{queuevacancy}, we obtain the following system of ordinary differential equations for link $a$:
\bsq\label{cltm-2}
\bqn
\der{\la_a(t)}t&=& \cas{{ll} k_a(L_a-V_a t,0) V_a -g_a(t), & t\leq \frac {L_a}{V_a} \\f_a(t-\frac{L_a}{V_a})-g_a(t), & t>\frac{L_a}{V_a}}\\
\der{\gamma_a(t)}t&=& \cas{{ll} -k_a(W_a t,0) W_a+K_a W_a -f_a(t), & t\leq \frac{L_a}{W_a}\\ g_a(t-\frac{L_a}{W_a})-f_a(t), &t>\frac{L_a}{W_a}}
\eqn
\esq
where the link demand and supply are given by \refe{def:ds}, and $f_a(t)$ and $g_a(t)$ can be computed by \refe{junctionmodel}. 
In particular, when link $a$ is initially empty, we have
\bsq \label{cltm-2empty}
\bqn
\der{\la_a(t)}t&=&\cas{{ll}0, & t\leq \frac{L_a}{V_a}\\ f_a(t-\frac{L_a}{V_a})-g_a(t), &t>\frac{L_a}{V_a}}\\
\der{\gamma_a(t)}t&=&\cas{{ll} K_aW_a-f_a(t), &t\leq \frac{L_a}{W_a}\\ g_a(t-\frac{L_a}{W_a})-f_a(t), & t>\frac{L_a}{W_a}}
\eqn
\esq
where the link demand and supply are given by 
\bqs
d_a(t)&=&\cas{{ll}0, & t\leq \frac{L_a}{V_a}\\\min\left\{ f_a(t-\frac{L_a}{V_a}) +H(\la_a(t)), C_a \right\},&t>\frac{L_a}{V_a}}\\
s_a(t)&=&\cas{{ll}C_a,&t\leq\frac{L_a}{W_a}\\\min\left\{ g_a(t-\frac{L_a}{W_a}) +H(\gamma_a(t)), C_a \right\}, &t>\frac{L_a}{W_a}}
\eqs

Furthermore, we can track the proportions of commodity $ p $ at the upstream and downstream boundaries of link $a$, which are denoted by $\eta_{a, p }(t)$ and $\xi_{a, p }(t)$, respectively. From \refe{junctionmodel} we can calculate the upstream commodity proportions as
\bsq \label{commodityproportion}
\bqn
\eta_{a, p }(t)&=&\frac{f_{a, p }(t)}{f_a(t)}.
\eqn
Since all vehicles follow the first-in-first-out principle in the road network, the downstream commodity proportions can be obtained from the upstream commodity proportions:
\bqn
\xi_{a, p }(t)&=&\eta_{a, p }(t- \pi _a(t)),
\eqn 
where $ \pi _a(t)$ is the travel time of vehicles arriving at the downstream boundary of link $a$ at time $t$:
\bqn
G_a(t)&=&F_a(t- \pi _a(t)).
\eqn 
In addition, the turning proportions can be calculated from the downstream commodity proportions:
\bqn
\xi_{a\to b}(t)&=&\sum_{ p \in  P _a \cap  P _b} \xi_{a, p }(t).
\eqn
\esq

Apparently both formulations, \refe{cltm-1} and \refe{cltm-2}, are systems of ordinary differential equations with delays, since link demands and supplies depend on historical states. Also in both formulations, the boundary fluxes $f_a(t)$ and $g_a(t)$ are important variables, which can be determined by the state variables. Once the cumulative flows $F_a(t)$ and $G_a(t)$ are found, traffic conditions inside a link can be obtained with Newell's model, \refe{u-domain-sol1}. In this sense, both Newell's model and LTM are complementary to each other in solving the traffic dynamics in a road network when origin demands $d_i(t)$ ($i\in I$), origin commodity proportions $\xi_{i,\o}(t)$ ($i\in I$), and destination supplies $s_o(t)$ ($o\in O)$ are all given.

\section{Stationary states in a network}
Under general boundary conditions at origins and destinations and initial conditions on links, LTM, either \refe{cltm-1} or \refe{cltm-2}, has to be numerically solved. In this section we are interested in analytically defining and solving the stationary patterns in LTM.
\subsection{Definition and properties of stationary states}
Traffic on a homogeneous link $a$ is considered stationary if for $x_a\in [0, L_a]$ (hereafter the subscript $a$ is omitted) 
\bqn
\pd{k(x,t)}t&=&0. \label{def:ss}
\eqn
In \citep{jin2012statics}, it was shown that in stationary states the flux is constant; i.e., for $x\in[0, L]$ and $t$
\bqn
q(x,t)=q.
\eqn
In addition, there are four types of stationary states: strictly under-critical (SUC), critical (C), strictly over-critical (SOC), and a zero-speed shock wave (ZS). 

\begin{theorem} A link can be stationary at C, SUC, SOC, or ZS states, in which the cumulative flow can be written as ($\beta\in [0,1]$)
\bsq \label{cfss}
\bqn
A(x,t)&=&N(x)+qt=qt+\cas{{ll} N(0)-k_1 x, & x\in[0,(1-\beta) L]\\N(L)+ (L-x) k_2, & x\in ((1-\beta) L, L]}
\eqn
where $k_1$ and $k_2$ ($k_1\leq \bar{K}\leq k_2$) are the respectively under- and over-critical densities corresponding to the flux $q$:
\bqn
q&=&k_1 V=(K-k_2) W. \label{def:k1k2}
\eqn
\esq
Here $\beta$ is the congested portion of a road.
\end{theorem}
{\em Proof}. From \refe{def:ss} and \refe{lwr_vt}, we can see that $\frac{\partial^2 A(x,t)}{\partial x \partial t}=0$ and $A_t-Q(-A_x)=0$, which lead to $\pd{A_t}x=0$, and $\pd{A_t}t=-Q_k(-A_x)\pd{A_x}t =0$. Thus $A_t=q$ is constant, and $A(x,t)=N(x)+qt$. From \refe{lwr_vt} we have $q-Q(-N_x)=0$. Thus $N_x(x)=k_1$ or $k_2$, which are given in \refe{def:k1k2}, and $N(x)$ is a piecewise linear function in $x$ with slopes of $-k_1$ and $-k_2$.
\ben
\item When $q=C$, we have $k_1=k_2=\bar{K}$, and traffic is stationary at a critical state with $N(x)=N(0)-\bar{K} x$, which corresponds to the C stationary state. 
\item When $q<C$, assuming that at a discontinuous point $y$, $N_x(x)=-k_2$ for $x<y$ and $N_x(x)=-k_1$ for $x>y$. Then for a small $t>0$ we have from \refe{hl-triangular}
\bqs
A(y,t)&=&N(y)+qt=\min\{N(y)+Ct, N(y-V t), N(y+W t)+KWt\}\\&=&\min\{N(y)+Ct, N(y)+k_2 Vt, N(y)-k_1 Wt +KWt\}.
\eqs
However, since $q<C$, we have that $Ct>qt$, $k_2Vt>qt$, and $(K-k_1)Wt>qt$. Thus the above equation does not have a solution. Therefore, it is impossible that the upstream part is more congested than the downstream part. 
In contrast, if at a discontinuous point $y$, $N_x(x)=-k_1$ for $x<y$ and $N_x(x)=-k_2$ for $x>y$. Then for a small $t>0$ from \refe{hl-triangular} we have
\bqs
A(y,t)&=&N(y)+qt=\min\{N(y)+Ct, N(y-V t), N(y+W t)+KWt\}\\&=&\min\{N(y)+Ct, N(y)+k_1 Vt, N(y)-k_2 Wt +KWt\}\\
&=&\min\{N(y)+Ct, N(y)+q t, N(y)+q t\}.
\eqs
Thus it is possible that the downstream part is more congested. Therefore, we can have three types of solutions when $q<C$: (i) $N(x)=N(0)-k_1 x$ ($x\in [0,L]$) for the SUC stationary state; (ii) $N(x)=N(0)-k_2 x$ ($x\in [0,L]$) for the SOC stationary state; or (iii) $N(x)=N(0)-k_1 x$ ($x\in [0,(1-\beta) L]$) and $N(x)=N(L)+ (L-x)k_2$ ($x\in ((1-\beta) L, L]$ and $\beta\in (0,1)$) for the ZS stationary state.  
\een

Furthermore, if we extend the range of $\beta$ to include $0$ and $1$, all of the four types of stationary states can be written as in \refe{cfss}.
In particular, when $q=C$, \refe{cfss} leads to the critical stationary state; when $q<C$ and $\beta=0$, \refe{cfss} leads to the SUC stationary state; when $q<C$ and $\beta=1$, \refe{cfss} leads to the SOC stationary state; and when $q<C$ and $\beta\in (0,1)$, \refe{cfss} leads to the ZS stationary state.
\eop

More properties of stationary states are presented in the following.
\begin{theorem} If link $a$ is stationary, then the in-flux and out-flux are equal and time-independent
\bqn
f(t)&=&g(t)=q, \label{fluxss}
\eqn
and both queue and vacancy sizes are time-independent
\bsq \label{qvss}
\bqn
\la(t)&=&\la= \beta  (1-\frac q C) KL,\\
\gamma(t)&=&\gamma=(1-\beta) (1-\frac q C) KL.
\eqn
\esq
\end{theorem}
{\em Proof}. In a stationary state, obviously \refe{fluxss} is true.
From \refe{queuevacancy} we have
\bqs
\la(t)&=&F(t-\frac LV)-G(t)=N(0)- k_1 L-N(L)\equiv \la,\\
\gamma(t)&=&G(t-\frac LW) +KL-F(t)=N(L)+k_2 L-N(0) \equiv \gamma,
\eqs
which lead to \refe{qvss}, since $k_2-k_1=(1-\frac qC)K$. In particular, we have the following cases: in C stationary states, $\la=\gamma=0$; in SUC stationary states, $\la=0$, and $\gamma=(1-\frac qC)K L>0$; in SOC stationary states, $\la=(1-\frac qC)K L>0$, and $\gamma=0$; in ZS stationary states, $\la= \beta  (1-\frac qC)K>0$, and $\gamma=(1-\beta) (1-\frac qC)K L>0$.
\eop

Note that, however,  \refe{fluxss} or \refe{qvss} are not sufficient conditions for stationary states. Consider a link connecting origin $i$ to destination $o$, and $d_i=s_o=q<C$. If the initial density is given by 
\bqs
k(x,0)&=&\cas{{ll} k_1, & x\in [0,L_1]\\ k_3, & x\in (L_1,L_2)\\ k_2, & x\in [L_2,L]}
\eqs
where $0<L_1<L_2<L$, $q=k_1 V=(K-k_2)W$, and $k_1 L\leq k_1 L_1+k_3 (L_2-L_1)+k_2 (L-L_2) \leq k_2 L$. In this case, $f(t)=g(t)=q$, but there can be shock and rarefaction waves inside the link. Thus \refe{fluxss} and \refe{qvss} cannot be used to determine stationary states.

\begin{corollary} In stationary states, the link demand and supply can be determined by $q$ and $(1-\beta)$:
\bsq \label{stationaryds}
\bqn
d&=&\min\{q+ H( \beta  (1-\frac{q}{C}) KL), C\},\\
s&=&\min\{q+ H((1-\beta)  (1-\frac{q}{C}) KL),C\}.
\eqn
\esq
\end{corollary}
These can be easily derived from  \refe{def:ds} and \refe{qvss}.

\subsection{Stationary states in a diverge-merge network}

\begin{figure}
\bc\includegraphics[height=1in]{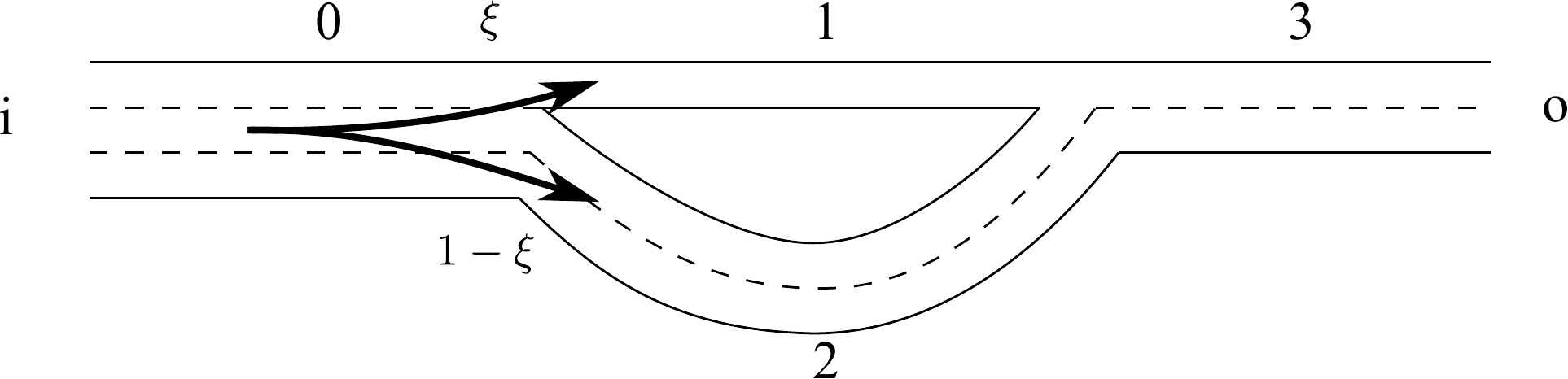}\ec
\caption{A diverge-merge network with one O-D pair and two intermediate links}\label{dm2network}
\end{figure}

 In this subsection, we try to find the stationary states in the diverge-merge network shown in \reff{dm2network}, where the capacities of four links are $(C_0, C_1, C_2, C_3)$, the turning proportion to link $1$ is $\xi$, the origin demand $d_i=C_0$, and the destination supply $s_o=C_3$. This is the so-called traffic statics problem in \citep{jin2012statics}, but we solve it with LTM.

In stationary states, link $0$ cannot be at SUC or ZS states; otherwise, $s_0=C_0$, and $f_0=\min\{d_i,s_0\}=C_0>q_0$. Therefore, link $0$ is at either SOC or C states, in which $\la_0>0$, and $d_0=C_0$.
Similarly, in stationary states, link $3$ is at either SUC or C states, in which $\gamma_3>0$, and $s_3=C_3$. Therefore, at the diverge junction, we have from \refe{junctionmodel}
\bsq \label{dmss}
\bqn
\theta_1&=&\min\{1, \frac{s_1}{C_0 \xi}, \frac{s_2}{(1-\xi) C_0 }\},\nonumber\\
g_0&=&C_0 \theta_1,\nonumber\\
f_1&=&\xi g_0=\min\{\xi C_0, s_1, \mu s_2\},\\
f_2&=&(1-\xi) g_0=\min\{(1-\xi) C_0,\frac 1\mu s_1, s_2\},
\eqn
where $\mu=\frac{\xi}{1-\xi}$;
at the merge junction, we have from \refe{junctionmodel}
\bqn
\theta_2&=&\min\{\max\{\frac{d_1}{C_1},\frac{d_2}{C_2}\},\max\{\frac{C_3}{C_1+C_2},\frac{C_3-d_2}{C_1},\frac{C_3-d_1}{C_2}\}\},\nonumber\\
g_1&=&\min\{d_1, C_1 \theta_2\}=\min\{d_1, \max \{ C_3-d_2, \nu C_3\}\},\\
g_2&=&\min\{d_2, C_2 \theta_2\}=\min\{d_2, \max\{ C_3-d_1, (1-\nu) C_3\}\},
\eqn
where $\nu=\frac{C_1}{C_1+C_2}$.
\esq
Further from \refe{fluxss}, we can replace $f_a$ and $g_a$ by $q_a$ ($a=1,2$) in \refe{dmss} and substitute $d_a$ and $s_a$ by \refe{stationaryds}. Then we obtain four equations with four unknown variables: $q_a$ and $\beta_a$ ($a=1,2$), which can be solved for given $C_0$, $C_1$, $C_2$, $C_3$, and $\xi$. We can follow \citep{jin2012statics} to show that the solutions exist but may not be unique.

Here we only consider the SOC-SUC stationary states on the two intermediate links. In this case, $\beta_1=1$, and $\beta_2=0$. From \refe{stationaryds} we have $d_1=C_1$, $s_1=q_1<C_1$, $d_2=q_2<C_2$, and $s_2=C_2$. From \refe{dmss} we have
\bqs
q_1&=&\min\{\xi C_0, q_1, \mu C_2\},\\
q_2&=&\min\{(1-\xi) C_0,\frac 1\mu q_1, C_2\},\\
q_1&=&\min\{C_1, \max \{ C_3-q_2, \nu C_3\}\},\\
q_2&=&\min\{q_2, \max\{ C_3-C_1, (1-\nu) C_3\}\},
\eqs
From the first equation we have $q_1\leq \min\{\xi C_0, \mu C_2\}$, and from the second equation we have $q_2=\frac 1\mu q_1$. From the third equation we have $q_1=\max\{C_3-q_2,\nu C_3\}<C_1$, which leads to $q_2>C_3-C_1$ and $\nu C_3< C_1$. From the fourth equation we have $q_2 \leq \max\{C_3-C_1, (1-\nu) C_3\}$, which leads to $q_2 \leq (1-\nu) C_3$ and $C_3-q_2 \geq \nu C_3$. Thus $q_1=C_3-q_2$.  Therefore the necessary conditions for the existence of SOC-SUC stationary states are that  $C_3<C_1+C_2$, $C_3\leq C_0$, and $1-\frac{C_2}{C_3}<\xi<\frac{C_1}{C_3}$, and the corresponding fluxs are $q_1=\xi C_3$ and $q_2=(1-\xi)C_3$. It can also be verified that the sufficient conditions for the existence of SOC-SUC stationary states in \refe{dmss} are
\bqn
C_3<C_1+C_2, C_3<C_0, 1-\frac{C_2}{C_3}<\xi<\frac{C_1}{C_3}. \label{socsuc}
\eqn
An example is for $(C_0,C_1,C_2,C_3)=(3,1,2,2)$ and $\xi\in(0,\frac 12)$, which was first studied in \citep[][Chapter 7]{jin2003dissertation}.

\section{Stability of stationary states}
In this section we consider the stability of the SOC-SUC stationary states in a diverge-merge network when \refe{socsuc} is satisfied. We apply small perturbations to both links 1 and 2, but the traffic conditions are still SOC and SUC, respectively. Therefore, we have $\la_1(t)>0$, $\gamma_1(t)=0$, $\la_2(t)=0$, and $\gamma_2(t)>0$. From \refe{queuevacancy} we can see that $F_1(t)=G_1(t-\frac {L_1}{W_1})+K_1 L_1$ and $G_2(t)=F_2(t-\frac {L_2}{V_2})$, which lead to 
\bsq \label{delayeqn}
\bqn
f_1(t)&=&g_1(t-\frac {L_1}{W_1}),\\
g_2(t)&=&f_2(t-\frac {L_2}{V_2}).
\eqn
\esq
Further from \refe{def:ds} we have $d_1(t)=C_1$, $s_1(t)=f_1(t)=g_1(t-\frac {L_1}{W_1})$, $d_2(t)=g_2(t)=f_2(t-\frac {L_2}{V_2})$, and $s_2(t)=C_2$.

Then from the junction model \refe{junctionmodel} we have the following relationships at the diverge:
\bqs
f_1(t)&=&\min\{\xi C_0, f_1(t), \mu C_2\},\\
f_2(t)&=&\min\{(1-\xi) C_0, \frac 1\mu f_1(t), C_2\},
\eqs
which lead to 
\bsq \label{dmeqn}
\bqn
f_2(t)&=&\frac 1\mu f_1(t).
\eqn
At the merge, we have
\bqs
g_1(t)&=&\min\{C_1, \max\{C_3-g_2(t),\nu C_3\}\},\\
g_2(t)&=&\min\{g_2(t),\max\{C_3-C_1, (1-\nu)C_3\}\},
\eqs
which lead to
\bqn
g_1(t)&=&C_3-g_2(t).
\eqn
\esq

Combining \refe{delayeqn} for the links and \refe{dmeqn} for the junctions, we then obtain the following equation:
\bqn
f_1(t)&=&C_3-\frac 1\mu f_1(t-T), \label{pm}
\eqn 
where $T=\frac {L_1}{W_1}+\frac {L_2}{V_2}$. This is equivalent to the Poincar\'e map in \citep{jin2013stability}, which was derived based on the circular information propagation in the network.
But here the Poincar\'e map, \refe{pm}, is directly derived from LTM. 

For \refe{pm}, the equilibrium point is $f_1(t)=\xi C_3$, which is the stationary flux on link 1. If we denote $\tilde f_1(t)=f_1(t)-\xi C_3$, then we have $\tilde f_1(t)= -\frac 1\mu \tilde f_1(t)$, whose equilibrium point is at 0. By analyzing the property of this map, we can then determine the stability of LTM.
\begin{theorem}
For a diverge-merge network satisfying \refe{socsuc}, the SOC-SUC stationary states and the equilibrium points of the Poincar\'e map, \refe{pm}, are stable when $\xi>\frac 12$, and unstable when $\xi\leq \frac 12$.
\end{theorem}
This is the same as Theorem 4.1 in \citep{jin2013stability}. Furthermore, we can also study the stability property of the stationary states with LTM in more general networks as in \citep{jin2013stability}.

\section{On non-invariant junction models}
In this section, we will solve LTM under the empty initial conditions and constant origin demands and destination supply in a merge network shown in \reff{merge-loading}. The demands of the two origins are denoted by $d_1^-$ and $d_2^-$, respectively, and the supply of the destination is denoted by $s_3^+$. Here we assume that $L_1=L_2=L_3=1$, $C_1=C_2=C_3=1$, $d_1^-=1$, $d_2^-=\frac 14$, and $s_3^+=1$. We will solve the traffic statics problem in this simple network with both invariant and non-invariant merge models.

\bfg\bc
\includegraphics[width=3in]{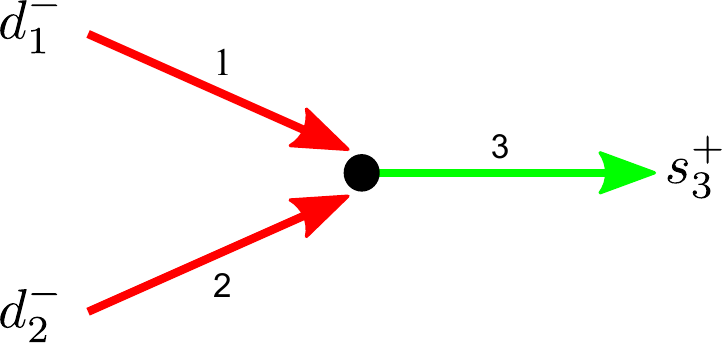}\caption{A merge network with constant origin demands and destination supply}\label{merge-loading}
\ec\efg

In stationary states, from \refe{stationaryds} we have $d_a=\min\{q_a+H(\beta_a(1-\frac {q_a}{C_a})K_a), C_a\}$ and $s_a=\min\{q_a+H((1-\beta_a)(1-\frac {q_a}{C_a})K_a),C_a\}$ ($a=1,2,3$). Here $q_a\leq C_a=1$. At the destination, $q_3=\min\{d_3,s_3^+\}=\min\{q_3+H(\beta_3(1-q_3)K_3), 1\}$, which leads to $\beta_3(1-q_3)=0$. Thus link 3 is stationary at C or SUC states: $q_3\leq 1$, and $s_3=1$. At the two origins, we have
\bqs
q_1&=&\min\{d_1^-,s_1\}=\min\{1,q_1+H((1-\beta_1)(1-q_1)K_1)\},\\
q_2&=&\min\{d_2^-,s_2\}=\min\{\frac 14,q_2+H((1-\beta_2)(1-q_2)K_2)\}.
\eqs
From the first equation, we have $(1-\beta_1)(1-q_1)=0$. Thus link 1 is stationary at C or SOC states: $q_1\leq 1$, and $d_1=1$. 
From the second equation, we can see that link 2 can be stationary at (1) SUC or ZS states: $\beta_2<1$, $q_2=\frac 14$, and 
\bsq \label{merge-link}
\bqn
d_2&=&\min\{\frac 14+H(\frac 34 \beta_2K_2),1\}=\cas{{ll} \frac 14, & \beta_2=0 \: (\m{SUC})\\1, &\beta_2\in(0,1) \: (\m{ZS})}
\eqn 
 or (2) SOC states: $\beta_2=1$, $q_2\leq \frac 14$, and 
\bqn
d_2&=&\min\{q_2+H((1-q_2)K_2),1\}=1.
\eqn 
\esq

First, we apply the invariant fair merge model in \refe{junctionmodel}:
\bsq\label{merge-junction1}
\bqn
q_1&=&\min\{d_1, \max \{ 1-d_2, \frac 12\}\}=\min\{1,\max\{1-d_2,\frac12\}\},\\
q_2&=&\min\{d_2, \max\{ 1-d_1, \frac 12\}\}=\min\{d_2,\frac 12\}.
\eqn
\esq
It can be verified that the solutions to \refe{merge-link} and \refe{merge-junction1} exist: $q_1=\frac 34$, $\beta_1=1$, $d_1=1$, $s_1=\frac 34$, $q_2=\frac 14$, $\beta_2=0$, $d_2=\frac 14$, and $s_2=1$. That is, link 1 is stationary at an SOC state, and link 2 at an SUC state. 

In contrast, we then apply the non-invariant fair merge model in \citep{jin2003merge}:
\bsq\label{merge-junction2}
\bqn
q_1&=&\frac{d_1}{d_1+d_2} \min\{d_1+d_2,s_3\}=\frac{1}{1+d_2},\\
q_2&=&\frac{d_2}{d_1+d_2} \min\{d_1+d_2,s_3\}=\frac{d_2}{1+d_2},
\eqn
\esq
since $d_1=s_3=1$. From \refe{merge-link}, there are three possible stationary states on link 2: (1) $\beta_2=0$, $d_2=\frac 14$, and $q_2=\frac 14$; (2) $\beta_2\in(0,1)$, $d_2=1$, and $q_2=\frac 14$; (3) $\beta_2=1$, $d_2=1$, and $q_2\leq 14$. However, neither of these satisfy (\ref{merge-junction2}b).
 Therefore, the traffic statics problem does not have a solution with the non-invariant merge model.

Note that, in \citep{jin2010merge}, it was shown that \refe{merge-junction1} and \refe{merge-junction2} are equivalent in continuous CTM, even though different interior states can arise around the merge. However, this example demonstrates that, for LTM to be well defined, we require junction models be invariant. Furthermore we make the following conjecture regarding the relationship between invariant junction models and well-defined traffic statics problem. 

\begin{conjecture} \label{conjecture}
A junction model is invariant if and only if solutions exist for the corresponding traffic statics problem of LTM for a junction network.
\end{conjecture}

\section{Conclusion}
In this paper, we first reviewed the Hamilton-Jacobi formulation of the LWR model and the corresponding Hopf-Lax formula, which consistent with the traditional variational principle, and derived Newell's simplified kinematic wave model inside a U-shaped spatial-temporal domain, where the boundary cumulative flows are given. We then applied the Hopf-Lax formula to define link demand and supply functions and used invariant junction models to calculate boundary fluxes. We also defined link queue and vacancy sizes and presented two continuous formulations of the link transmission model (LTM).
We further defined stationary states on a link and solved the stationary states in a diverge-merge network with constant origin demands, destination supplies, and commodity proportions. We applied LTM to directly derive a Poincar\'e map to analyze the stability of stationary states in a diverge-merge network. Furthermore we compared invariant and non-invariant merge models and showed that LTM is not well-defined with non-invariant junction models.

From this study, we can see that Newell's model is based on the traditional variational principle or Hopf-Lax formula, which can be used to solve the Cauchy-Dirichlet problem with given boundary cumulative flows. In contrast, LTM can be used to solve the Cauchy-Neumann problem with given boundary demand and supply functions and macroscopic junction models. These two models are complementary to each other: we can use LTM to find the boundary cumulative flows of a link and then use Newell's model to find traffic conditions inside the link. Thus  LTM, complemented by Newell's model, is equivalent to the network kinematic wave model, which is a system of partial differential equations.

LTM can be extended for other bottlenecks and traffic systems, including capacity drop, heterogeneous traffic, and signalized intersections. LTM with general initial conditions can be useful for studying traffic dynamics in closed networks. In addition, the computational efficiency and accuracy of LTM also warrant more studies. 

In addition, more properties of LTM can be investigated. For examples, we will be interested in proving Conjecture \ref{conjecture}, which, if true, can serve as another definition of invariant junction models. In addition, we will be interested in proving the existence, studying the stability, and developing algorithms of stationary states in large-scale networks with LTM.

Finally, we will also be interested in applying LTM to design and analyze traffic control strategies and  solve the dynamic traffic assignment problem analytically and numerically.

\section*{Acknowledgments} We would like to thank Zhe Sun and Dr. Yifeng Yu of UC Irvine and Dr. Ke Han of Imperial College London for discussions.


\begin{thebibliography}{38}
\expandafter\ifx\csname natexlab\endcsname\relax\def\natexlab#1{#1}\fi
\expandafter\ifx\csname url\endcsname\relax
  \def\url#1{\texttt{#1}}\fi
\expandafter\ifx\csname urlprefix\endcsname\relax\def\urlprefix{URL }\fi

\bibitem[{Bardos et~al.(1979)Bardos, Leroux, and Nedelec}]{bardos1979first}
Bardos, C., Leroux, A., Nedelec, J., 1979. {First order quasilinear equations
  with boundary conditions}. Communications in partial differential equations
  4~(9), 1017--1034.

\bibitem[{Beckmann(1952)}]{beckmann1952transportation}
Beckmann, M., 1952. A continuous model of transportation. Econometrica: Journal
  of the Econometric Society 20~(4), 643--660.

\bibitem[{Claudel and Bayen(2010)}]{claudel2010lax2}
Claudel, C.~G., Bayen, A.~M., 2010. {Lax--Hopf based incorporation of internal
  boundary conditions into Hamilton-Jacobi equation. part II: Computational
  methods}. IEEE Transactions on Automatic Control 55~(5), 1158--1174.

\bibitem[{Daganzo(1995)}]{daganzo1995ctm}
Daganzo, C.~F., 1995. {The cell transmission model {II}: Network traffic}.
  Transportation Research Part B 29~(2), 79--93.

\bibitem[{Daganzo(2005{\natexlab{a}})}]{daganzo2005variationalKW}
Daganzo, C.~F., 2005{\natexlab{a}}. A variational formulation of kinematic
  waves: basic theory and complex boundary conditions. Transportation Research
  Part B 39~(2), 187--196.

\bibitem[{Daganzo(2005{\natexlab{b}})}]{daganzo2005variationalKW2}
Daganzo, C.~F., 2005{\natexlab{b}}. A variational formulation of kinematic
  waves: Solution methods. Transportation Research Part B 39~(10), 934--950.

\bibitem[{Daganzo(2006)}]{daganzo2006variational}
Daganzo, C.~F., 2006. {On the variational theory of traffic flow:
  well-posedness, duality and applications}. Networks and Heterogeneous Media
  1~(4), 601--619.

\bibitem[{Daganzo and Geroliminis(2008)}]{daganzo2008analytical}
Daganzo, C.~F., Geroliminis, N., 2008. An analytical approximation for the
  macroscopic fundamental diagram of urban traffic. Transportation Research
  Part B 42~(9), 771--781.

\bibitem[{Evans(1998)}]{evans1998pde}
Evans, L., 1998. {Partial Differential Equations}. American Mathematical
  Society.

\bibitem[{Fl{\"o}tter{\"o}d and Rohde(2011)}]{flotterod2011operational}
Fl{\"o}tter{\"o}d, G., Rohde, J., 2011. Operational macroscopic modeling of
  complex urban road intersections. Transportation Research Part B 45~(6),
  903--922.

\bibitem[{Garavello and Piccoli(2006)}]{garavello2006tfn}
Garavello, M., Piccoli, B., 2006. {Traffic Flow on Networks}. Vol.~1. Applied
  Mathematics Series.

\bibitem[{Gentile et~al.(2007)Gentile, Meschini, and
  Papola}]{gentile2007macroscopic}
Gentile, G., Meschini, L., Papola, N., 2007. {Spillback congestion in dynamic
  traffic assignment: A macroscopic flow model with time-varying bottlenecks}.
  Transportation Research Part B 41~(10), 1114--1138.

\bibitem[{Han et~al.(2012)Han, Piccoli, Friesz, and Yao}]{han2012continuous}
Han, K., Piccoli, B., Friesz, T., Yao, T., 2012. A continuous-time link-based
  kinematic wave model for dynamic traffic networks. arXiv preprint
  arXiv:1208.5141.

\bibitem[{Hidas(2005)}]{hidas2005modelling}
Hidas, P., 2005. {Modelling vehicle interactions in microscopic simulation of
  merging and weaving}. Transportation Research Part C 13~(1), 37--62.

\bibitem[{Ho and Wong(2006)}]{ho2006continuum}
Ho, H., Wong, S., 2006. Two-dimensional continuum modeling approach to
  transportation problems. Journal of Transportation Systems Engineering and
  Information Technology 6~(6), 53--68.

\bibitem[{Holden and Risebro(1995)}]{holden1995unidirection}
Holden, H., Risebro, N.~H., 1995. A mathematical model of traffic flow on a
  network of unidirectional roads. SIAM Journal on Mathematical Analysis
  26~(4), 999--1017.

\bibitem[{Jin(2003)}]{jin2003dissertation}
Jin, W.-L., 2003. Kinematic wave models of network vehicular traffic. Ph.D.
  thesis, University of California, Davis.
\newline\urlprefix\url{http://arxiv.org/abs/math.DS/0309060}

\bibitem[{Jin(2010)}]{jin2010merge}
Jin, W.-L., 2010. {Continuous kinematic wave models of merging traffic flow}.
  Transportation Research Part B 44~(8-9), 1084--1103.

\bibitem[{Jin(2012{\natexlab{a}})}]{jin2012_riemann}
Jin, W.-L., 2012{\natexlab{a}}. {A Riemann solver for a system of hyperbolic
  conservation laws at a general road junction}. Arxiv preprint.
\newline\urlprefix\url{http://arxiv.org/abs/1204.6727}

\bibitem[{Jin(2012{\natexlab{b}})}]{jin2012network}
Jin, W.-L., 2012{\natexlab{b}}. A kinematic wave theory of multi-commodity
  network traffic flow. Transportation Research Part B 46~(8), 1000--1022.

\bibitem[{Jin(2012{\natexlab{c}})}]{jin2012kinematic}
Jin, W.-L., 2012{\natexlab{c}}. A kinematic wave theory of multi-commodity
  network traffic flow. Transportation Research Part B 46~(8), 1000--1022.

\bibitem[{Jin(2012{\natexlab{d}})}]{jin2012_link}
Jin, W.-L., 2012{\natexlab{d}}. A link queue model of network traffic flow.
  arXiv preprint arXiv:1209.2361.

\bibitem[{Jin(2012{\natexlab{e}})}]{jin2012statics}
Jin, W.-L., 2012{\natexlab{e}}. The traffic statics problem in a road network.
  Transportation Research Part B 46~(10), 1360--1373.

\bibitem[{Jin(2013)}]{jin2013stability}
Jin, W.-L., 2013. {Stability and bifurcation in network traffic flow: A
  Poincar\'e map approach}. Transportation Research Part B 57, 191--208.

\bibitem[{Jin and Zhang(2003)}]{jin2003merge}
Jin, W.-L., Zhang, H.~M., 2003. On the distribution schemes for determining
  flows through a merge. Transportation Research Part B 37~(6), 521--540.

\bibitem[{Laval et~al.(2012)Laval, He, and Castrillon}]{laval2012stochastic}
Laval, J.~A., He, Z., Castrillon, F., 2012. Stochastic extension of newell's
  three-detector method. Transportation Research Record: Journal of the
  Transportation Research Board 2315, 73--80.

\bibitem[{Lebacque(1996)}]{lebacque1996godunov}
Lebacque, J.~P., 1996. {The Godunov scheme and what it means for first order
  traffic flow models}. Proceedings of the 13th International Symposium on
  Transportation and Traffic Theory, 647--678.

\bibitem[{Lebacque and Khoshyaran(2005)}]{lebacque2005network}
Lebacque, J.~P., Khoshyaran, M., 2005. {First order macroscopic traffic flow
  models: Intersection modeling, Network modeling}. Proceedings of the 16th
  International Symposium on Transportation and Traffic Theory, 365--386.

\bibitem[{LeVeque(2002)}]{leveque2002fvm}
LeVeque, R.~J., 2002. Finite volume methods for hyperbolic problems. Cambridge
  University Press, Cambridge; New York.

\bibitem[{Lighthill and Whitham(1955)}]{lighthill1955lwr}
Lighthill, M.~J., Whitham, G.~B., 1955. {On kinematic waves: II. A theory of
  traffic flow on long crowded roads}. Proceedings of the Royal Society of
  London A 229~(1178), 317--345.

\bibitem[{Moskowitz(1965)}]{moskowitz1965discussion}
Moskowitz, K., 1965. { Discussion of `freeway level of service as in uenced by
  volume and capacity characteristics' by D.R. Drew and C. J. Keese}. Highway
  Research Record 99, 43--44.

\bibitem[{Newell(1993)}]{newell1993sim}
Newell, G.~F., 1993. {A simplified theory of kinematic waves in highway traffic
  {I}: General theory. {II}: Queuing at freeway bottlenecks. {III}:
  Multi-destination flows}. Transportation Research Part B 27~(4), 281--313.

\bibitem[{Nie and Zhang(2002)}]{nie2002sqm}
Nie, X., Zhang, H., 2002. {The Formulation of A Link Based Dynamic Network
  Loading Model Considering Queue Spillovers}. Tech. rep., working Paper
  UCD-ITS-Zhang-2002-6.

\bibitem[{Richards(1956)}]{richards1956lwr}
Richards, P.~I., 1956. Shock waves on the highway. Operations Research 4~(1),
  42--51.

\bibitem[{Tamp{\`e}re et~al.(2011)Tamp{\`e}re, Corthout, Cattrysse, and
  Immers}]{tampere2011generic}
Tamp{\`e}re, C., Corthout, R., Cattrysse, D., Immers, L., 2011. A generic class
  of first order node models for dynamic macroscopic simulation of traffic
  flows. Transportation Research Part B 45~(1), 289--309.

\bibitem[{Yperman(2007)}]{yperman2007link}
Yperman, I., 2007. The link transmission model for dynamic network loading.
  Ph.D. thesis.

\bibitem[{Yperman et~al.(2006)Yperman, Logghe, Tampere, and
  Immers}]{yperman2006mcl}
Yperman, I., Logghe, S., Tampere, C., Immers, B., 2006. {The Multi-Commodity
  Link Transmission Model for Dynamic Network Loading}. Proceedings of the TRB
  Annual Meeting.

\bibitem[{Zhang et~al.(2013)Zhang, Nie, and Qian}]{zhang2013modelling}
Zhang, H., Nie, Y., Qian, Z., 2013. Modelling network flow with and without
  link interactions: the cases of point queue, spatial queue and cell
  transmission model. Transportmetrica B: Transport Dynamics 1~(1), 33--51.

\end{thebibliography}
\end{document}